\documentclass[review,onefignum,onetabnum]{siamart250211}

\usepackage{booktabs}
\usepackage{amsfonts}
\usepackage{graphicx}
\usepackage{epstopdf}
\usepackage{algorithmic}
\usepackage{mathdots}
\usepackage{multirow}
\usepackage[shortlabels]{enumitem}
\usepackage{comment}
\usepackage{rotating}
\usepackage{xcolor}
\usepackage{titlecaps}
\usepackage{caption}
\usepackage{subcaption}

\newsiamremark{remark}{Remark}
\newsiamremark{example}{Example}

\newsiamremark{hypothesis}{Hypothesis}
\crefname{hypothesis}{Hypothesis}{Hypotheses}
\newsiamthm{claim}{Claim}

\newcommand{\mat}[1]{\begin{bmatrix}#1 \\ \end{bmatrix}}

\newcommand{\correc}[1]{{\color{black}  #1}}

\usepackage{titlecaps}

\headers{Preconditioner for parabolic optimal control problems}{S. Y. Hon, P. Y. Fung, and X.-L. Lin}

\title{An optimal \correc{diagonalization-based} preconditioner for parabolic optimal control problems}

\author{Sean Y. Hon\footnotemark[2] \and Po Yin Fung\footnotemark[3] \and  Xue-lei Lin\footnotemark[1]
	\footnotetext{\footnotemark[1]~\lowercase{\MakeUppercase Corresponding \MakeUppercase Author. \MakeUppercase School of \MakeUppercase Science, \MakeUppercase Harbin \MakeUppercase Institute of \MakeUppercase Technology, \MakeUppercase Shenzhen 518055, \MakeUppercase  China (\email{linxuelei@hit.edu.cn}).}}
	\footnotetext{\footnotemark[2]~\lowercase{\MakeUppercase Department of \MakeUppercase Mathematics, \MakeUppercase Hong \MakeUppercase Kong \MakeUppercase Baptist \MakeUppercase University, \MakeUppercase Kowloon \MakeUppercase Tong, \MakeUppercase Hong \MakeUppercase Kong \MakeUppercase{SAR} (\email{seanyshon@hkbu.edu.hk}).}}
	\footnotetext{\footnotemark[3]~\lowercase{\MakeUppercase Department of \MakeUppercase Mathematics, \MakeUppercase Hong \MakeUppercase Kong \MakeUppercase Baptist \MakeUppercase University, \MakeUppercase Kowloon \MakeUppercase Tong, \MakeUppercase Hong \MakeUppercase Kong \MakeUppercase{SAR} (\email{pyfung@hkbu.edu.hk}).}}
}

\usepackage{amsopn}

\begin{document}
	\nolinenumbers
	\maketitle
	
	\begin{abstract}
		In this work, we propose a novel diagonalization-based preconditioner for the all-at-once linear system arising from the optimal control problem of parabolic equations. The proposed preconditioner is constructed based on an $\epsilon$-circulant modification to the rotated block diagonal (RBD) preconditioning technique \correc{and can be} efficiently diagonalized by fast Fourier transforms in a parallel-in-time fashion. To our knowledge, this marks the first application of the $\epsilon$-circulant modification to RBD preconditioning. Before our work, the studies of \correc{parallel-in-time} preconditioning techniques for the optimal control problem are mainly focused on $\epsilon$-circulant modification to Schur complement based preconditioners, which involves multiplication of forward and backward evolutionary processes and thus square the condition number. Compared with those Schur complement based preconditioning techniques in the literature, the advantage of the proposed $\epsilon$-circulant modified RBD preconditioning is that it does not involve the multiplication of forward and backward evolutionary processes. When the generalized minimal residual method is deployed on the preconditioned system, we prove that when choosing $\epsilon=\mathcal{O}(\sqrt{\tau})$ with $\tau$ being the temporal step-size\correc{, }the convergence rate of the preconditioned GMRES solver is independent of the matrix size and the regularization parameter. \correc{Numerical results are provided to demonstrate the effectiveness of our proposed solvers.}
	\end{abstract}
	
	
	\section{Introduction}\label{sec:introduction}
	Developing diagonalization-based parallel-in-time (PinT) \correc{methods} for solving optimization problems constrained by partial differential equations (PDEs) has gained substantial attention in recently years. We refer to \cite{lions_1971, Hinze_2008, Troltzsch2010, Borzi_2011} and the references therein for a comprehensive overview of these optimization problems. Various efficient PinT preconditioners have been proposed to solve all-at-once linear systems arising from time-dependent PDEs \cite{doi:10.1137/16M1062016,McDonald2017,doi:10.1137/20M1316354,doi:10.1137/19M1309869,LiLinHon2023,vcaklovic2023parallel} and PDE-constrained optimal control problems \cite{LevequePearson22,liuWu_optimal,WuWangZhou2023,WuZhou2020,hondongSC2023,Bouillon2024,FungHon_2024}; \correc{for a review paper on these diagonalization-based methods, see \cite{2020arXiv200509158G} and the references therein reported.}
	
	In this work, we consider solving the distributed optimal control problem constrained by a parabolic equation. Namely, the following quadratic cost functional is minimized:
	\begin{equation}\label{eqn:Cost_functional_heat}
		\min_{y,u}~ \mathcal{J}(y,u):=\frac{1}{2}\| y - g \|^{2}_{L^2(\Omega \times (0,T))} + \frac{\gamma}{2}\| u \|^{2}_{L^2(\Omega \times (0,T))},
	\end{equation}
	constrained by a parabolic equation with certain initial and boundary conditions
	\begin{equation}\label{eqn:heat}
		\left\{
		\begin{array}{lc}
			y_{t} - \mathcal{L} y = f + u, \quad (x,t)\in \Omega \times (0,T], \qquad  y = 0, \quad (x,t)\in \partial \Omega \times (0,T], \\
			y(x,0)=y_0, \quad x \in \Omega,
		\end{array}
		\right.\,
	\end{equation}
	where $u,g \in \correc{L^2(\Omega \times (0,T))}$ are the distributed control and the targeted tracking trajectory, respectively, $\gamma>0$ is a regularization parameter, $\mathcal{L}=\nabla \cdot (a({x})\nabla )$, and both $f$ and $y_0$ are given functions. The theoretical aspects of existence, uniqueness, and regularity of the solution, under suitable assumptions, were thoroughly studied in \cite{lions_1971}. \correc{Following a standard argument \cite{lions_1971}, the} optimal solution of \correc{(\ref{eqn:Cost_functional_heat}) - (\ref{eqn:heat})} can be characterized by the following system:
	\begin{equation}\label{eqn:heat_2}
		\left\{
		\begin{array}{lc}
			y_{t} - \mathcal{L} y - \frac{1}{\gamma} p= f,\quad (x,t)\in \Omega \times (0,T], \qquad y = 0,\quad  (x,t)\in \partial \Omega \times (0,T], \\
			y(x,0)=y_0, \quad  x \in \Omega,\\
			-p_{t} - \mathcal{L} p + y = g, \quad (x,t)\in \Omega \times (0,T], \qquad p = 0, \quad (x,t) \in \partial \Omega \times (0,T], \\
			p(x,T) = 0, \quad  x \in \Omega,
		\end{array}
		\right.\,
	\end{equation}
	where \correc{$p$ denotes the Lagrange multiplier, also known as the co-state variable. In deriving the second PDE in \eqref{eqn:heat_2}, the control variable $u$ is eliminated using the optimality condition $\gamma u - p = 0$. For further details, we refer the reader to \cite[Section 3.1]{WuWangZhou2023} and \cite{lions_1971}.}

	\correc{Following \cite{LevequePearson22}, we discretize (\ref{eqn:heat_2}) using the backward Euler method in time with a uniform step size $\tau=T/N$, and apply some spatial discretization with mesh size parameter $m$, which gives
	\begin{eqnarray*}
		M_{m}\frac{\mathbf{y}_m^{(k_1+1)} - \mathbf{y}_m^{(k_1)} }{\tau} +  K_m\mathbf{y}_m^{(k_1+1)}  &=&  M_{m}\left(\mathbf{f}_m^{(k_1+1)} + \frac{1}{\gamma} \mathbf{p}_m^{(k_1+1)}\right), \qquad k_1=0,1,2,\dots, N-1,\\
		-M_{m}\frac{\mathbf{p}_m^{(k_2+1)} - \mathbf{p}_m^{(k_2)} }{\tau} +  K_m \mathbf{p}_m^{(k_2)}    &=&  M_{m} \left( \mathbf{g}_m^{(k_2)} - \mathbf{y}_m^{(k_2)} \right),  \qquad k_2=0,1,2,\dots N-1,
	\end{eqnarray*} where the matrices $M_{m}$ and $K_{m}$ represent the mass matrix and the stiffness matrix, respectively, if a finite element method is used. For the finite difference method, we have $M_{m}=I_{m}$ and $K_{m}=-L_{m}$, where $-L_{m}$ is the discretization matrix of the negative Laplacian.
		
	Combining the given initial and final conditions that $\mathbf{y}_m^{(0)}=y_0$ and $\mathbf{p}_m^{(N)}=\mathbf{0}$, one needs to solve the following all-at-once system assembled using the formulas above for $k_1 = 0, 1, \dots, N - 2$ and $k_2 = 1, 2, \dots, N - 1$ (see, e.g., \cite{Bouillon2024, WuZhou2020}, where the system is derived in the same way):
	\begin{equation}\label{eqn:main_system_before}
		\widetilde{\mathcal{A}} \begin{bmatrix} \mathbf{y}\\ \mathbf{p} \end{bmatrix} = \begin{bmatrix} \mathbf{g}\\  \mathbf{f} \end{bmatrix},
	\end{equation}
	where we have $\mathbf{y} = [ \mathbf{y}_m^{(1)}, \cdots, \mathbf{y}_m^{(N-1)}]^{\top}$, $\mathbf{p} = [ \mathbf{p}_m^{(1)},\cdots, \mathbf{p}_m^{(N-1)}]^{\top}$, 
	\begin{equation*}
		\mathbf{f}=\mat{ \tau M_m\mathbf{f}_m^{(1)} + M_{m}\mathbf{y}_m^{(0)}\\
			\tau M_m \mathbf{f}_m^{(2)}\\
			\vdots\\
			\tau M_m \mathbf{f}_m^{(N-1)}}, \quad
		\mathbf{g} = \mat{\tau M_m \mathbf{g}_m^{(1)}\\
			\vdots\\
			\tau M_m \mathbf{g}_m^{(N-2)}\\
			\tau M_m \mathbf{g}_m^{(N-1)}+M_m\mathbf{p}_m^{(N)}},
	\end{equation*}
		\begin{eqnarray}\label{eqn:matrix_A_before}
			\widetilde{\mathcal{A}} &=&
			\begin{bmatrix} 
				\tau I_{N-1} \otimes M_{m}  &   B_{N-1}^{\top} \otimes M_{m} + \tau  I_{N-1} \otimes K_m \\
				B_{N-1} \otimes M_{m} + \tau  I_{N-1} \otimes K_m  &  -\frac{\tau}{\gamma} I_{N-1} \otimes M_{m} 
			\end{bmatrix},
		\end{eqnarray}
		and the matrix $B_{N-1} \in \mathbb{R}^{N-1 \times N-1}$ is as follows.
		\begin{equation}\label{eqn:Toeplitz_Mat_B}
			B_{N-1} = \begin{bmatrix}
				1 &   &  &  & \\
				-1  & 1    & & & \\
				&  -1  & 1  & &  \\
				&     & \ddots & \ddots &  \\
				&  &    & -1 & 1
			\end{bmatrix}.
		\end{equation} After solving the system \eqref{eqn:main_system_before} for $\mathbf{y}$ and $\mathbf{p}$, one proceeds to compute $\mathbf{y}_m^{(N)}$ and $\mathbf{p}_m^{(0)}$ via the following discretization formulas for $k_1=N-1$ and $k_2=0$:
        \begin{eqnarray*}
		(M_{m}+\tau K_m) \mathbf{y}_m^{(N)} &=& M_m \mathbf{y}_m^{(N-1)} + \tau M_m \mathbf{f}_m^{(N)}+\frac{\tau}{\gamma}M_m\mathbf{p}_m^{(N)},\\
	(M_{m}+\tau K_m) \mathbf{p}_m^{(0)} &=& M_m \mathbf{p}_m^{(1)} + \tau M_m \mathbf{g}_m^{(0)}-{\tau}M_m\mathbf{y}_m^{(0)}.
	\end{eqnarray*}
        
        The two linear systems above are classical elliptic problems, which can be efficiently solved with linear complexity by solvers proposed in the literature, see, e.g., \cite{bornemann1996,fulton1986multigrid} and the references therein.}
		
        
        \correc{Throughout}, the following conditions on both $M_{m}$ and $K_{m}$ are assumed, which are in line with \cite[Section 2]{doi:10.1137/20M1316354}:
		\begin{enumerate}
			\item The matrices $M_{m}$ and $K_{m}$ are (real) symmetric positive definite;
			\item The matrices $M_{m}$ and $K_{m}$ are sparse, having only $\mathcal{O}(m)$ nonzero elements;
			\item The matrix $M_{m}$ has a condition number uniformly upper bounded  with respect to $m$, i.e., there exists a positive finite constant $c_0$ independent of $m$ and $n$ such that 
			\begin{equation}\label{eqn:constant_c0}
				\sup\limits_{m>0}\kappa_2(M_{m})\leq c_0<+\infty,
			\end{equation}
			where $\kappa_2(M_{m}):=||M_{m}||_2||M_m^{-1}||_2$ denotes the condition number of $M_m$. 
		\end{enumerate}
		In fact, these assumptions are easily met by applying typical finite element discretization to the spatial term $\mathcal{L}$. 
		
		\correc{By letting $n=N-1$, we} further transform (\ref{eqn:matrix_A_before}) into the following equivalent system\correc{, for which our proposed preconditioner is specifically developed:}
		\begin{equation}\label{eqn:main_system}
			{\mathcal{A}}\underbrace{\begin{bmatrix} \sqrt{\gamma} \mathbf{y}\\ \mathbf{p} \end{bmatrix}}_{\correc{=:}{\bf x}} = \underbrace{\begin{bmatrix} \mathbf{g}\\ -\sqrt{\gamma} \mathbf{f} \end{bmatrix}}_{\correc{=:}{\bf b}},
		\end{equation}
		where 
		\begin{eqnarray}\nonumber
			{\mathcal{A}} &=&
			\begin{bmatrix} 
				\alpha {I}_n \otimes M_m   & B_n^{\top}\otimes M_m + \tau  {I}_n \otimes K_{m}\\
				-(B_n\otimes M_m + \tau  {I}_n\otimes K_{m})  &  \alpha {I}_n \otimes M_m 
			\end{bmatrix}\\\label{eqn:matrix_A}
			&=&
			\begin{bmatrix} 
				\alpha {I}_n \otimes M_m & \mathcal{T}^{\top}\\
				-\mathcal{T}  & \alpha {I}_n \otimes M_m 
			\end{bmatrix}.
		\end{eqnarray}
		Note that {$\alpha = \frac{\tau}{\sqrt{\gamma}}$}, $I_k$ denotes the $k \times k$ identity matrix, and
		\begin{equation}\label{eqn:matrix_T}
			\mathcal{T}=B_n\otimes M_m + \tau  {I}_n\otimes K_{m}.
		\end{equation}

					In what follows, we will develop \correc{a preconditioner for the nonsymmetric system (\ref{eqn:main_system}), to be applied within GMRES \correc{\cite{Saad_Schultz_1986}}.}
					
					For $\mathcal{A}$, we first consider the following rotated block diagonal (RBD) preconditioner \correc{firstly} proposed in \cite{bailu2021}:
					\begin{eqnarray}\label{eqn:matrix_P}
						\mathcal{P} = \frac{1}{2}
						\begin{bmatrix} 
							\mathcal{T}^{\top} + \alpha {I}_n \otimes M_m  & \\
							& \mathcal{T} + \alpha {I}_n \otimes M_m 
						\end{bmatrix}
						\begin{bmatrix} 
							I_{mn}  &  I_{mn} \\
							-I_{mn}  &  I_{mn} 
						\end{bmatrix},
					\end{eqnarray}
					where $\mathcal{T}$ is defined in (\ref{eqn:matrix_T}). As will be shown in Theorem \ref{thm:spectum_of_H_inv_B}, the matrix $\mathcal{P} $ can achieve an excellent preconditioning effect for $\mathcal{A}$. Yet, $\mathcal{P}$ is in general computationally expensive to invert \correc{in a PinT manner}. Thus, we propose the following PinT preconditioner \correc{(denoted} by RBD $\epsilon$-circulant preconditioner) as a modification of $\mathcal{P}$, which can be efficiently implemented \correc{in parallel}:

					\begin{eqnarray}\label{eqn:matrix_P_eps}
						\mathcal{P}_{\epsilon} = \frac{1}{2}
						\begin{bmatrix} 
							\mathcal{C}^{\top}_{\epsilon}+ \alpha {I}_n \otimes M_m  & \\
							& \mathcal{C}_{\epsilon} + \alpha {I}_n \otimes M_m 
						\end{bmatrix}
						\begin{bmatrix} 
							I_{mn}  &  I_{mn} \\
							-I_{mn}  &  I_{mn} 
						\end{bmatrix},
					\end{eqnarray}
					where 
					\begin{eqnarray*}\label{eqn:matrix_C_eps}
						\mathcal{C}_{\epsilon} = C_{\epsilon,n} \otimes M_m + \tau  {I}_n\otimes K_{m}.
					\end{eqnarray*}
					Note that $C_{\epsilon,n}$ is defined as 
					\begin{eqnarray}\label{eqn:matrix_C_eps_n}
						C_{\epsilon,n}=\begin{bmatrix}
							1 &   &  &  & -\epsilon\\
							-1  & 1    & & & \\
							&  -1  & 1  & &  \\
							&     & \ddots & \ddots &  \\
							&  &    & -1 & 1
						\end{bmatrix},
					\end{eqnarray}
					and $\epsilon \in (0,1]$. Since $C_{\epsilon,n}$ is an $\epsilon$-circulant matrix, it admits the following decomposition \cite{BiniLatoucheMeini,Bertaccini_Ng_2003}
					\begin{equation}\label{eqn:decomposition_eps}
						C_{\epsilon,n} = D_{\epsilon}^{-1} \mathbb{F}_{n} \Lambda_{\epsilon,n} \mathbb{F}_{n}^{*} D_{\epsilon},
					\end{equation}
					where \[D_{\epsilon}={\rm diag} \left(\epsilon^{\frac{i-1}{\correc{n}}}\right)_{i=1}^{n},\] \[\mathbb{F}_n=\frac{1}{\sqrt{n}}[\theta_n^{(i-1)(j-1)}]_{i,j=1}^{n} \quad \textrm{with} \quad  \theta_n =\exp\left(\frac{2\pi {\bf i}}{n}\right),\] and \[\Lambda_{\epsilon,n}= {\rm diag}(\lambda_{i-1}^{(\epsilon)})_{i=1}^{n}\quad \textrm{with} \quad\lambda_{k}^{(\epsilon)}=1-\epsilon^{\frac{1}{n}} \theta_{n}^{-k}.\]
					
					Several efficient preconditioned iterative solvers have been proposed for \eqref{eqn:matrix_A_before}. One notable class of such methods is the circulant-based preconditioners \cite{WuZhou2020,Bouillon2024,FungHon_2024}, whose main idea is to replace the Toeplitz matrix $B_n$ defined by \eqref{eqn:Toeplitz_Mat_B} by a circulant or skew-circulant matrix in order to form an efficient preconditioner. Despite their numerically observed success, theoretically these preconditioners have not shown to be parameter-robust \correc{ in theory, meaning that their performance can degrade significantly as problem parameters—such as the mesh size or the regularization parameter, primarily} due to the low-rank \correc{term in} the circulant-based approximation. 
					
					Another major existing method is preconditioning technique based on Schur complement approximations \correc{\cite{pearson2012regularization,WuWangZhou2023,LevequePearson22,Linheatopt2022}}, which is a typical preconditioning approach in the context of preconditioning for saddle point systems (see, e.g., \correc{\cite{MurphyGolubWathen2000,ANU:9672992})} and their effectiveness is based on  approximating the Schur complements by multiplications of easily invertible matrices (e.g., \correc{\cite{PearsonWathen2012,LevequePearson22,AxelssonNeytcheva2006}}). Recently, $\epsilon$-circulant modifications have been introduced to  such \textcolor{black}{approximate-Schur-complement} preconditioning techniques for the optimal control problems, which leads to a significant improvement on computational efficiency and parallelism along \correc{the} time dimension (see, e.g, \cite{doi:10.1137/20M1316354, doi:10.1137/19M1309869, LiLinHon2023, Linhonwave2022, WuWangZhou2023, 2020arXiv200509158G}).
					

                    \correc{Our proposed preconditioning method, while still based on $\epsilon$-circulant matrices, diverges from the approximate-Schur-complement approach. Specifically, the Toeplitz matrix $B_n$ is \correc{replaced with an $\epsilon$-circulant matrix, as detailed in \eqref{eqn:matrix_C_eps_n}}. A key advancement of our method is that it achieves an optimal convergence rate when $\epsilon = \mathcal{O}(\sqrt{\tau})$. This sets it apart from the $\epsilon$-circulant modified-approximate-Schur-complement preconditioners proposed in~\cite{Linheatopt2022, WuWangZhou2023}.
                    
                    We note that both~\cite{Linheatopt2022, WuWangZhou2023} focus on the second-order Crank--Nicolson scheme, while our work concerns the backward Euler scheme. On one hand, it remains unclear whether their methods can be directly applied to our setting. On the other hand, although the time discretizations differ, the underlying challenge of constructing effective $\epsilon$-circulant-based preconditioners is closely related. In particular, the choice of $\epsilon$ plays a critical role in balancing approximation quality and numerical stability in both contexts. Therefore, comparing the required conditions on $\epsilon$ remains meaningful. Specifically, the method in~\cite{Linheatopt2022} requires a more restrictive condition $\epsilon = \mathcal{O}(\tau^2)$. While the approach in~\cite{WuWangZhou2023} also allows $\epsilon = \mathcal{O}(\sqrt{\tau})$, it imposes a strong theoretical constraint that $n$ must be even—a limitation that our method does not impose.
                    }

                    \correc{Furthermore, it is worth noting that our preconditioning approach differs from that in \cite{LevequePearson22}, even though both resulting linear systems are derived using the same discretization scheme. Specifically, by enforcing initial and final time conditions as $(M_m+\tau K_m)\mathbf{y}_m^{(0)}=(M_m+\tau K_m){y}_0$ and $(M_m+\tau K_m)\mathbf{y}_m^{(N)}=\mathbf{0}$, the system in \cite{LevequePearson22} is formulated as:
                    \begin{equation}\label{eqn:LevequePearson_system}
                            \begin{bmatrix}
                            \tau \hat{I}_{N+1} \otimes M_m &  B_{N+1}^\top \otimes M_m + \tau I_{N+1}\otimes K_m \\ B_{N+1} \otimes M_m + \tau I_{N+1}\otimes K_m & -\frac{\Delta t}{\gamma} \check{I}_{N+1} \otimes M_m \\
                        \end{bmatrix}
                        \begin{bmatrix}
                            \boldsymbol{y}\\ \boldsymbol{p}  \\
                        \end{bmatrix}
                        =
                        \begin{bmatrix}
                            \boldsymbol{g}\\ \boldsymbol{f}  \\
                        \end{bmatrix},
                    \end{equation}
                    with $\boldsymbol{y} = [ \mathbf{y}_m^{(0)}, \mathbf{y}_m^{(1)}, \cdots, \mathbf{y}_m^{(N)}]^{\top}$, $\boldsymbol{p} = [ \mathbf{p}_m^{(0)},\cdots, \mathbf{p}_m^{(N)}]^{\top}$, 
                	\begin{equation*}
                		\boldsymbol{f}=\mat{ (M_m+\tau K_m)y_0\\
                			\tau M_m \mathbf{f}_m^{(1)}\\
                			\vdots\\
                			\tau M_m \mathbf{f}_m^{(N-1)}}, \quad
                		\boldsymbol{g} = \mat{\tau M_m \mathbf{g}_m^{(0)}\\
                			\vdots\\
                			\tau M_m \mathbf{g}_m^{(N-1)}\\
                			 \mathbf{0}},
                	\end{equation*} and
                    \[
                    \hat{I}_{N+1}
                    =
                    \begin{bmatrix}
                           1 & & & \\ & \ddots & &\\ & & 1 & \\& &   & 0 \\
                    \end{bmatrix} \in \mathbb{R}^{N+1 \times N+1},
                    \quad\quad
                    \check{I}_{N+1}
                    =
                    \begin{bmatrix}
                           0 & & & \\ & 1 & &\\ & & \ddots & \\& &   & 1 \\
                    \end{bmatrix}\in \mathbb{R}^{N+1 \times N+1}.
                    \]
                    A Schur-complement-based preconditioner was then proposed in \cite{LevequePearson22} for the system \eqref{eqn:LevequePearson_system}; however, its parallelization was not discussed. In contrast, our proposed preconditioner for the system \eqref{eqn:main_system_before} can be implemented efficiently using fast Fourier transforms and is well-suited for parallel execution.
                    
                    }

					\begin{remark}
						In addition to the RBD preconditioner of interest in this work, other effective preconditioners such as the preconditioned square block matrix \cite{Axelsson2020} could potentially be incorporated with the $\epsilon$-circulant preconditioning technique. We have chosen the RBD preconditioner due to its simple form, whose inversion \correc{is} only mainly determined by that of the block diagonal with $\mathcal{C}_{\epsilon}+ \alpha {I}_n \otimes M_m$ and $\mathcal{C}^{\top}_{\epsilon}+ \alpha {I}_n \otimes M_m$. Such an inversion can be straightforwardly implemented in a PinT way.
					\end{remark}
					
					The paper is organized as follows. In Section \ref{sec:main}, we provide our main results on the spectral analysis for our proposed preconditioners. Numerical examples are given in Section \ref{sec:numerical} for supporting the performance of our proposed preconditioner. As last, conclusions are given in Section \ref{sec:conclusion}.

					\section{Main Result}\label{sec:main}
					Before presenting our main preconditioning result, we introduce some useful notation in what follows. 
					
					Let \[\mathcal{G} = \frac{1}{2}\begin{bmatrix} 
						I_{mn}  &  I_{mn} \\
						-I_{mn}  &  I_{mn} 
					\end{bmatrix},\qquad \mathcal{H}_{\epsilon}=\begin{bmatrix} 
						\mathcal{C}_{\epsilon}^{\top} + \alpha {I}_n \otimes M_m  & \\
						& \mathcal{C}_{\epsilon} + \alpha {I}_n \otimes M_m 
					\end{bmatrix},\]~\textcolor{black}{and \[
						\mathcal{B}=\begin{bmatrix} 
							\mathcal{T}^{\top} + \alpha {I}_n \otimes M_m  & \mathcal{T}^{\top} - \alpha {I}_n \otimes M_m\\
							-\mathcal{T} + \alpha {I}_n \otimes M_m  & \mathcal{T} + \alpha {I}_n \otimes M_m 
						\end{bmatrix}.
						\]
						Then, it is straightforward to verify that
						\begin{equation*}
							\mathcal{A}=\mathcal{B}\mathcal{G},\qquad \mathcal{P}_{\epsilon}=\mathcal{H}_{\epsilon}\mathcal{G}.
						\end{equation*}
						To check the invertibility of $\mathcal{P}_{\epsilon}$, it suffices to verify the invertibility of $\mathcal{G}$ and $\mathcal{H}_{\epsilon}$ respectively. 
						Clearly, $\mathcal{G}$ is invertible. By \eqref{eqn:decomposition_eps}, it is straightforward to see that $\mathcal{C}_{\epsilon}+\alpha I_n\otimes M_m$ is similar to a block diagonal matrix with each diagonal block having the form $\alpha M_m + \correc{ \tau  K_{m}}+ z_1 I_m+z_2{\bf i}I_m$ ($z_1\geq 0$, $z_2\in\mathbb{R}$). Consequently, $\mathcal{C}_{\epsilon}+\alpha I_n\otimes M_m$ is invertible. Moreover, $\mathcal{C}_{\epsilon}^{\top} + \alpha {I}_n \otimes M_m$ as the transpose of $\mathcal{C}_{\epsilon} + \alpha {I}_n \otimes M_m$ is also invertible. Thus, $\mathcal{H}_{\epsilon}$, as a block diagonal matrix with $\mathcal{C}_{\epsilon}^{\top} + \alpha {I}_n \otimes M_m$ and $\mathcal{C}_{\epsilon} + \alpha {I}_n \otimes M_m$ as its diagonal blocks, is also invertible. Therefore, $\mathcal{P}_{\epsilon}$ is invertible.}  
					
					In this work, we propose the use of the matrix $\mathcal{P}_{\epsilon}$ to precondition the saddle point system described in \eqref{eqn:main_system}. In other words, we shall apply the GMRES solver to solve the following preconditioned system
					\begin{equation}\label{eqn:preconditioned_system}
						\mathcal{P}_{\epsilon}^{-1}\mathcal{A}{\bf x}=\mathcal{P}_{\epsilon}^{-1}{\bf b}.
					\end{equation}
					We will focus on investigating the convergence rate of the GMRES solver for the preconditioned system mentioned above, by theoretically estimating a suitable choice of $\epsilon$.
					
					Denote \correc{with}
					\begin{align*}
						& \mathcal{W} = \begin{bmatrix}
							I_{n} \otimes M_{m}& \\
							& I_{n} \otimes M_{m}
						\end{bmatrix} \qquad  \quad \widetilde{\mathcal{T}}=B_n\otimes I_m + \tau  {I}_n\otimes (M_{m}^{-\frac{1}{2}} K_{m} M_{m}^{-\frac{1}{2}}),\\
						&\widetilde{\mathcal{B}}=\begin{bmatrix} 
							\widetilde{\mathcal{T}}^{\top} + \alpha {I}_{mn}   & \widetilde{\mathcal{T}}^{\top} - \alpha {I}_{mn}\\
							-\widetilde{\mathcal{T}} + \alpha {I}_{mn}  &\widetilde{\mathcal{T}} + \alpha {I}_{mn} 
						\end{bmatrix},\qquad\widetilde{\mathcal{C}}_{\epsilon}=C_{\epsilon,n}\otimes I_m + \tau  {I}_n\otimes (M_{m}^{-\frac{1}{2}} K_{m} M_{m}^{-\frac{1}{2}}),\\
						&\quad \widetilde{\mathcal{H}}_{\epsilon}=\begin{bmatrix} 
							\widetilde{\mathcal{C}}_{\epsilon}^{\top} + \alpha {I}_{mn} & \\
							& \widetilde{\mathcal{C}}_{\epsilon} + \alpha {I}_{mn} 
						\end{bmatrix},\qquad \widetilde{\mathcal{P}}_{\epsilon}=\widetilde{\mathcal{H}}_{\epsilon}\mathcal{G}.
					\end{align*}
					It is then easy to see that
					\begin{equation*}
						\mathcal{T}=\mathcal{W}^{\frac{1}{2}}\widetilde{\mathcal{T}}\mathcal{W}^{\frac{1}{2}},\quad \mathcal{B}=\mathcal{W}^{\frac{1}{2}}\widetilde{\mathcal{B}}\mathcal{W}^{\frac{1}{2}},\quad \mathcal{C}_{\epsilon}=\mathcal{W}^{\frac{1}{2}}\widetilde{\mathcal{C}}_{\epsilon}\mathcal{W}^{\frac{1}{2}},\quad\mathcal{H}_{\epsilon}=\mathcal{W}^{\frac{1}{2}}\widetilde{\mathcal{H}}_{\epsilon}\mathcal{W}^{\frac{1}{2}}.
					\end{equation*}
					With the equations above, it is straightforward to verify that \eqref{eqn:preconditioned_system} can be rewritten as \[\mathcal{G}^{-1}\mathcal{W}^{-\frac{1}{2}}\widetilde{\mathcal{H}}_{\epsilon}^{-1}\widetilde{\mathcal{B}}\mathcal{W}^{\frac{1}{2}}\mathcal{G}{\bf x}=\mathcal{G}^{-1}\mathcal{W}^{-\frac{1}{2}}\widetilde{\mathcal{H}}_{\epsilon}^{-1}\mathcal{W}^{-\frac{1}{2}}{\bf b},\]  meaning that \eqref{eqn:preconditioned_system} can be equivalently converted into the following system
					\begin{equation}\label{auxiliarysystem}
						\widetilde{\mathcal{H}}_{\epsilon}^{-1}\widetilde{\mathcal{B}}\tilde{\bf x}=  \tilde{\bf b},
					\end{equation}
					where  $\tilde{\bf b}=\widetilde{\mathcal{H}}_{\epsilon}^{-1}\mathcal{W}^{-\frac{1}{2}}{\bf b}$; the unknown $\tilde{\bf x}$ in \eqref{auxiliarysystem} and the unknown  ${\bf x}$ in \eqref{eqn:preconditioned_system} are related by
					\begin{equation*}
						\tilde{\bf x}=    \mathcal{W}^{\frac{1}{2}}\mathcal{G}{\bf x}.
					\end{equation*}
					
					We will examine the theoretical convergence rate of the GMRES solver for the linear system \eqref{auxiliarysystem} and subsequently relate it to the convergence rate of the GMRES solver for \eqref{eqn:preconditioned_system} using Theorem \ref{gmrescvgrelationthm}.
					
					The convergence behavior of \correc{the} GMRES solver is closely related to \correc{Krylov subspaces}. For a square matrix $\mathbf{E}\in \mathbb{R}^{k \times k}$ and a vector ${\bf y}\in \mathbb{R}^{k\times 1}$, \correc{the} Krylov subspace of degree $j\geq 1$ \correc{generated by the matrix $\mathbf{E}$ and the vector ${\bf y}$} is defined as follows
					$$\mathcal{K}_j(\mathbf{E},{\bf y}):=\text{span}\{{\bf y},\mathbf{E}{\bf y},\mathbf{E}^2{\bf y},\dots,\mathbf{E}^{j-1}{\bf y} \}.$$
					For a set $\mathcal{S}$ and a point $z$, we denote
					$$z+\mathcal{S}:=\{z+y|y\in \mathcal{S}\}.$$
					We recall the relation between the iterative solution by GMRES and the Krylov subspace in the following lemma.
					\begin{lemma}\correc{(see, e.g.,~\cite{saad2003iterative})}\label{krylov}
						For a non-singular $k\times k$ real linear system ${\bf Z} {\bf y} = {\bf w}$, let ${\bf y}_j$ be the GMRES iterative solution at the $j$-th iteration $(j \geq 1)$ with ${\bf y}_0$ as an initial guess. Then, the $j$-th iteration solution ${\bf y}_j$ \correc{minimizes} the residual error over the Krylov subspace $\mathcal{K}_j({\bf Z},{\bf r}_0)$ with ${\bf r}_0= {\bf w}-{\bf Z}{\bf y}_0,$ i.e.,	
						$$
						\mathbf{y}_{j}=\underset{\mathbf{c} \in \mathbf{y}_{0}+\mathcal{K}_{j}\left(\mathbf{Z}, \mathbf{r}_{0}\right)}{\arg \min }\|\mathbf{w}-\mathbf{Z} \mathbf{c}\|_{2}.
						$$
					\end{lemma}
					
					\correc{The following theorem establishes a connection between the convergence rate estimates of the GMRES solver applied to \eqref{auxiliarysystem} and \eqref{eqn:preconditioned_system}.}
                    
					\begin{theorem}\label{gmrescvgrelationthm}
						Let $\tilde{\mathbf{x}}_{0}$ be arbitrary initial guess for \eqref{auxiliarysystem}. Let $\mathbf{x}_{0}:=\mathcal{G}^{-1}\mathcal{W}^{-\frac{1}{2}} \tilde{\mathbf{x}}_{0}$ be the initial guess for \eqref{eqn:preconditioned_system}. Let $\mathbf{x}_{j}~(\tilde{\mathbf{x}}_{j},~{respectively})$ be the j-th $(j \geq 1)$ iteration solution derived by applying GMRES solver to \eqref{eqn:preconditioned_system} $( \eqref{auxiliarysystem}, ~{respectively})$ with $\correc{\mathbf{x}_{0}}~(\correc{\tilde{\mathbf{x}}_{0}},~{respectively})$ as an initial guess. Then,
						$$
						\left\|\mathbf{r}_{j}\right\|_{2} \leq \sqrt{2}||\mathcal{W}^{-\frac{1}{2}}||_2||\tilde{\bf r}_j||_2,\quad j=1,2,\dots,
						$$
						where $\mathbf{r}_{j}:=\mathcal{P}_{\epsilon}^{-1} \mathbf{b}-\mathcal{P}_{\epsilon}^{-1} \mathcal{A} \mathbf{x}_{j}~(\tilde{\mathbf{r}}_{j}:=\tilde{\bf b} - \widetilde{\mathcal{H}}_{\epsilon}^{-1}\widetilde{\mathcal{B}} \tilde{\mathbf{x}}_{j}$, respectively$)$ denotes the residual vector at the $j$-th GMRES iteration for \eqref{eqn:preconditioned_system} $($\eqref{auxiliarysystem}, respectively$)$.
					\end{theorem}
					\begin{proof}
						Note that
						\begin{equation*}
							\mathcal{G}^{-1}\mathcal{W}^{-\frac{1}{2}}\tilde{\bf x}_j-{\bf x}_0= \mathcal{G}^{-1}\mathcal{W}^{-\frac{1}{2}}\tilde{\bf x}_j-\mathcal{G}^{-1}\mathcal{W}^{-\frac{1}{2}}\tilde{\bf x}_0=\mathcal{G}^{-1}\mathcal{W}^{-\frac{1}{2}}(\tilde{\bf x}_j-\tilde{\bf x}_0).
						\end{equation*}
						By definition of $\tilde{\bf x}_j$ and $\tilde{\bf x}_0$, we know that $\tilde{\bf x}_j-\tilde{\bf x}_0\in\mathcal{K}_j(\widetilde{\mathcal{H}}_{\epsilon}^{-1}\widetilde{\mathcal{B}},\tilde{\mathbf{r}}_{0})$ with $\tilde{\mathbf{r}}_{0}=\tilde{\bf b} - \widetilde{\mathcal{H}}_{\epsilon}^{-1}\widetilde{\mathcal{B}} \tilde{\mathbf{x}}_{0}$, meaning that there exists real numbers $\eta_0,\eta_1,\dots,\eta_{j-1}$ such that
						\begin{align*}
							\tilde{\bf x}_j-\tilde{\bf x}_0&=   \sum\limits_{k=0}^{j-1}\eta_k(\widetilde{\mathcal{H}}_{\epsilon}^{-1}\widetilde{\mathcal{B}})^{k}\tilde{\mathbf{r}}_{0}\\
							&=   \sum\limits_{k=0}^{j-1}\eta_k\mathcal{W}^{\frac{1}{2}}\mathcal{G}(\mathcal{G}^{-1}\mathcal{W}^{-\frac{1}{2}}\widetilde{\mathcal{H}}_{\epsilon}^{-1}\widetilde{\mathcal{B}}\mathcal{W}^{\frac{1}{2}}\mathcal{G})^{k}\mathcal{G}^{-1}\mathcal{W}^{-\frac{1}{2}}\tilde{\mathbf{r}}_{0}\\
							&= \sum\limits_{k=0}^{j-1}\eta_k\mathcal{W}^{\frac{1}{2}}\mathcal{G}(\mathcal{P}_{\epsilon}^{-1} \mathcal{A})^{k}\mathcal{G}^{-1}\mathcal{W}^{-\frac{1}{2}}\tilde{\mathbf{r}}_{0}\\
							&= \sum\limits_{k=0}^{j-1}\eta_k\mathcal{W}^{\frac{1}{2}}\mathcal{G}(\mathcal{P}_{\epsilon}^{-1} \mathcal{A})^{k}(\mathcal{P}_{\epsilon}^{-1}{\bf b}-\mathcal{P}_{\epsilon}^{-1}\mathcal{A}{\bf x}_0)= \mathcal{W}^{\frac{1}{2}}\mathcal{G}\sum\limits_{k=0}^{j-1}\eta_k(\mathcal{P}_{\epsilon}^{-1} \mathcal{A})^{k}{\bf r}_0.
						\end{align*}
						Therefore, we have $\mathcal{G}^{-1}\mathcal{W}^{-\frac{1}{2}}(\tilde{\bf x}_j-\tilde{\bf x}_0)
						\in\mathcal{K}_j(\mathcal{P}_{\epsilon}^{-1} \mathcal{A},{\bf r}_0)$. Hence, 
						\begin{equation*}
							\mathcal{G}^{-1}\mathcal{W}^{-\frac{1}{2}}\tilde{\bf x}_j-{\bf x}_0\in \mathcal{K}_j(\mathcal{P}_{\epsilon}^{-1} \mathcal{A},{\bf r}_0).
						\end{equation*}
						In other words, $\mathcal{G}^{-1}\mathcal{W}^{-\frac{1}{2}}\tilde{\bf x}_j\in {\bf x}_0+ \mathcal{K}_j(\mathcal{P}_{\epsilon}^{-1} \mathcal{A},{\bf r}_0)$.
						Then, Lemma \ref{krylov} implies that
						\begin{align*}
							||\mathbf{r}_{j}||_2&=||\mathcal{P}_{\epsilon}^{-1} \mathbf{b}-\mathcal{P}_{\epsilon}^{-1} \mathcal{A} \mathbf{x}_{j}||_2 \\
							&\leq ||\mathcal{P}_{\epsilon}^{-1} \mathbf{b}-\mathcal{P}_{\epsilon}^{-1} \mathcal{A} \mathcal{G}^{-1}\mathcal{W}^{-\frac{1}{2}}\tilde{\bf x}_j||_2\\
							&=||\mathcal{G}^{-1}\mathcal{W}^{-\frac{1}{2}}(\widetilde{\mathcal{H}}_{\epsilon}^{-1}\mathcal{W}^{-\frac{1}{2}}{\bf b}-\widetilde{\mathcal{H}}_{\epsilon}^{-1}\widetilde{\mathcal{B}}\tilde{\bf x}_j)||_2\\
							&=||\mathcal{G}^{-1}\mathcal{W}^{-\frac{1}{2}}\tilde{\bf r}_j||_2\leq ||\mathcal{G}^{-1}||_2||\mathcal{W}^{-\frac{1}{2}}||_2||\tilde{\bf r}_j||_2=\sqrt{2}||\mathcal{W}^{-\frac{1}{2}}||_2||\tilde{\bf r}_j||_2.
						\end{align*}
						The proof is complete.
					\end{proof}
					
					
					According to Theorem \ref{gmrescvgrelationthm}, any convergence rate estimation for the GMRES solver applied to \eqref{auxiliarysystem} is also applicable to \eqref{eqn:preconditioned_system}. Namely, the GMRES convergence for \eqref{eqn:preconditioned_system} is guaranteed by its convergence for \eqref{auxiliarysystem}. As a result, in the following discussion, we will focus on investigating the convergence of the GMRES solver for \eqref{auxiliarysystem} by first examining the spectrum of its coefficient matrix, $\widetilde{\mathcal{H}}_{\epsilon}^{-1}\widetilde{\mathcal{B}}$.
					
					Denote \correc{with}
					\begin{equation*}
						\widetilde{\mathcal{H}}=\begin{bmatrix} 
							\widetilde{\mathcal{T}}^{\top} + \alpha {I}_{mn}& \\
							& \widetilde{\mathcal{T}} + \alpha {I}_{mn}
						\end{bmatrix}.
					\end{equation*}
					Clearly, $\widetilde{\mathcal{H}}_{\epsilon}$ is an approximation to $\widetilde{\mathcal{H}}$. In the next, we will study the spectrum of $\widetilde{\mathcal{H}}^{-1}\widetilde{\mathcal{B}}$ and then relate the spectrum of $\widetilde{\mathcal{H}}_{\epsilon}^{-1}\widetilde{\mathcal{B}}$  to the spectrum of $\widetilde{\mathcal{H}}^{-1}\widetilde{\mathcal{B}}$  by  the fact $\lim\limits_{\epsilon\rightarrow 0^{+}}\widetilde{\mathcal{H}}_{\epsilon}=\widetilde{\mathcal{H}}$.
					
					The following theorem shows that $ \widetilde{\mathcal{H}}$ is an ideal preconditioner for $ \widetilde{\mathcal{B}}$. We remark that the proof mirrors that of \cite[Theorem 4.1]{bailu2021}, with modifications. However, for self-containment, we decide to provide a complete proof below.
					
					Denote \correc{with} $\sigma({\bf C})$, the spectrum of a square matrix ${\bf C}$. 
					
					\begin{theorem}\label{thm:spectum_of_H_inv_B}
						The matrix $\widetilde{\mathcal{H}}^{-1}\widetilde{\mathcal{B}}$ is unitarily diagonalizable and its    eigenvalues satisfy $\sigma(\widetilde{\mathcal{H}}^{-1}\widetilde{\mathcal{B}}) \subseteq \{1+{\bf i}x|x\in[-1,1]\}$, where $\mathbf{i}=\sqrt{-1}$ is the imaginary unit.
					\end{theorem}
					\begin{proof}
						Consider
						\begin{eqnarray*}
							&&\widetilde{\mathcal{H}}^{-1}\widetilde{\mathcal{B}}\\
							&=& 
							\begin{bmatrix} 
								\widetilde{\mathcal{T}}^{\top} + \alpha {I}_{mn}  & \\
								& \widetilde{\mathcal{T}} + \alpha {I}_{mn}
							\end{bmatrix}^{-1}
							\begin{bmatrix} 
								\widetilde{\mathcal{T}}^{\top} + \alpha {I}_{mn}   & \widetilde{\mathcal{T}}^{\top} - \alpha {I}_{mn}\\
								-\widetilde{\mathcal{T}} + \alpha {I}_{mn}  &\widetilde{\mathcal{T}} + \alpha {I}_{mn} 
							\end{bmatrix}\\
							&=& 
							\begin{bmatrix}
								I_{mn} & (\widetilde{\mathcal{T}}^{\top}+ \alpha I_{mn})^{-1} (\widetilde{\mathcal{T}}^{\top} - \alpha I_{mn})\\
								(\widetilde{\mathcal{T}}+ \alpha I_{mn})^{-1} (-\widetilde{\mathcal{T}} + \alpha I_{mn}) & I_{mn}
							\end{bmatrix}\\
							&=&
							\begin{bmatrix}
								I_{mn} &  (\widetilde{\mathcal{T}}^{\top}-\alpha I_{mn})(\widetilde{\mathcal{T}}^{\top}+ \alpha I_{mn})^{-1} \\
								(\widetilde{\mathcal{T}}+ \alpha I_{mn})^{-1} (-\widetilde{\mathcal{T}} + \alpha I_{mn}) & I_{mn}
							\end{bmatrix},
						\end{eqnarray*}
						where the last equality comes from the following fact
						
						\begin{eqnarray*}
							&&(\widetilde{\mathcal{T}}^{\top}+ \alpha I_{mn})^{-1} (\widetilde{\mathcal{T}}^{\top} - \alpha I_{mn})\\
							&=& (\widetilde{\mathcal{T}}^{\top}+ \alpha I_{mn})^{-1} (\widetilde{\mathcal{T}}^{\top} - \alpha I_{mn})I_{mn}\\
							&=&(\widetilde{\mathcal{T}}^{\top}+ \alpha I_{mn})^{-1} (\widetilde{\mathcal{T}}^{\top} - \alpha I_{mn})[(\widetilde{\mathcal{T}}^{\top}+ \alpha I_{mn})(\widetilde{\mathcal{T}}^{\top}+ \alpha I_{mn})^{-1}]\\
							&=&(\widetilde{\mathcal{T}}^{\top}+ \alpha I_{mn})^{-1} [(\widetilde{\mathcal{T}}^{\top} - \alpha I_{mn})(\widetilde{\mathcal{T}}^{\top}+ \alpha I_{mn})](\widetilde{\mathcal{T}}^{\top}+ \alpha I_{mn})^{-1}\\
							&=&(\widetilde{\mathcal{T}}^{\top}+ \alpha I_{mn})^{-1} [(\widetilde{\mathcal{T}}^{\top}+ \alpha I_{mn})(\widetilde{\mathcal{T}}^{\top} - \alpha I_{mn})](\widetilde{\mathcal{T}}^{\top}+ \alpha I_{mn})^{-1}\\
							&=&[(\widetilde{\mathcal{T}}^{\top}+ \alpha I_{mn})^{-1} (\widetilde{\mathcal{T}}^{\top}+ \alpha I_{mn})](\widetilde{\mathcal{T}}^{\top} - \alpha I_{mn})(\widetilde{\mathcal{T}}^{\top}+ \alpha I_{mn})^{-1}\\
							&=&(\widetilde{\mathcal{T}}^{\top} - \alpha I_{mn})(\widetilde{\mathcal{T}}^{\top}+ \alpha I_{mn})^{-1}.
						\end{eqnarray*}
						Considering $\widetilde{\mathcal{E}} = (\widetilde{\mathcal{T}}+ \alpha I_{mn})^{-1} (-\widetilde{\mathcal{T}} + \alpha I_{mn})$ with its singular value decomposition, namely, $\widetilde{\mathcal{E}} = \mathcal{U} \Sigma \mathcal{V}^{\correc{\top}}$. Then, we can further decompose $\widetilde{\mathcal{H}}^{-1}\widetilde{\mathcal{B}}$ as follows:
						\begin{eqnarray*}
							\widetilde{\mathcal{H}}^{-1}\widetilde{\mathcal{B}}
							&=& 
							\begin{bmatrix}
								I_{mn}&-\widetilde{\mathcal{E}}^{\top}\\\widetilde{\mathcal{E}}&I_{mn}
							\end{bmatrix}= 
							\begin{bmatrix}
								\mathcal{V}& \\ &\mathcal{U}
							\end{bmatrix}
							\begin{bmatrix}
								I_{mn} & -\correc{\Sigma}\\\Sigma & I_{mn}
							\end{bmatrix}
							\begin{bmatrix}
								\mathcal{V}^{\correc{\top}}& \\ &\mathcal{U}^{\correc{\top}}
							\end{bmatrix}.
						\end{eqnarray*}
						Denote by 
						\begin{eqnarray*}
							\mathcal{Q}=\frac{1}{\sqrt{2}}\begin{bmatrix}
								I_{mn}&I_{mn}\\-\mathbf{i} I_{mn}& \mathbf{i} I_{mn}
							\end{bmatrix},
						\end{eqnarray*}
						where $\mathbf{i}=\sqrt{-1}$ is the imaginary unit. Note that $\mathcal{Q}$ is a unitary matrix. It holds that
						\begin{eqnarray*}
							\begin{bmatrix}
								I_{mn}&-\correc{\Sigma}\\\Sigma&I_{mn}
							\end{bmatrix}
							&=& \mathcal{Q}\begin{bmatrix}
								I_{mn}+\mathbf{i} \Sigma & \\ &I_{mn} - \mathbf{i} \Sigma
							\end{bmatrix}
							\mathcal{Q}^{*}.
						\end{eqnarray*}
						Denote $\mathcal{X} =
						\begin{bmatrix}
							\mathcal{V}& \\ &\mathcal{U}
						\end{bmatrix} \mathcal{Q}$.
						Consequently, it follows that
						\begin{eqnarray*}
							\widetilde{\mathcal{H}}^{-1}\widetilde{\mathcal{B}} &=& 
							\begin{bmatrix}
								\mathcal{V}& \\ &\mathcal{U}
							\end{bmatrix}
							\mathcal{Q}\begin{bmatrix}
								I_{mn}+\mathbf{i} \Sigma & \\ &I_{mn} - \mathbf{i} \Sigma
							\end{bmatrix}
							\mathcal{Q}^{*}
							\correc{\begin{bmatrix}
								\mathcal{V}^{\correc{\top}}& \\ &\mathcal{U}^{\correc{\top}}
							\end{bmatrix}}\\
							&=&\mathcal{X}\begin{bmatrix}
								I_{mn}+\mathbf{i} \Sigma & \\ &I_{mn} - \mathbf{i} \Sigma
							\end{bmatrix}\mathcal{X}^{*}
						\end{eqnarray*}
						Clearly, the eigenvalues of  $\widetilde{\mathcal{H}}^{-1}\widetilde{\mathcal{B}}$ locates in $\{1+{\bf i}x|x\in[-||\widetilde{\mathcal{E}}||_2,||\widetilde{\mathcal{E}}||_2]\}\subset\mathbb{C}$.
						As the matrix $\widetilde{\mathcal{T}}$ is positive definite (i.e., $\widetilde{\mathcal{T}} + \widetilde{\mathcal{T}}^\top$ is \correc{symmetric} positive definite), we know from \correc{\cite[Theorem~2.2]{bai2014}} that $\|\widetilde{\mathcal{E}}\|_{2} < 1$. \correc{Namely, $\widetilde{\mathcal{E}} =: V(\widetilde{\mathcal{T}}, \alpha)$, which is defined as the so-called extrapolated Cayley transform in \cite{bai2014}, is strictly less than one in the Euclidean norm. For more details, refer to the proof of \cite[Theorem~2.2]{bai2014}.} Therefore, $$\sigma(\widetilde{\mathcal{H}}^{-1}\widetilde{\mathcal{B}})\subseteq\{1+{\bf i}x|x\in[-||\widetilde{\mathcal{E}}||_2,||\widetilde{\mathcal{E}}||_2]\}\subseteq \{1+{\bf i}x|x\in[-1,1]\}.$$
						
						Moreover, it is trivial to see that $\mathcal{X}$ is unitary. Hence, $\widetilde{\mathcal{H}}^{-1}\widetilde{\mathcal{B}}$ is unitarily diagonalizable.
						The proof is complete.
					\end{proof}

					It is straightforward to verify that  $\widetilde{\mathcal{T}}+\alpha I_{mn}$ can be rewritten as 
					\begin{equation*}
						\widetilde{\mathcal{T}}+\alpha I_{mn}=\left[
						\begin{array}[c]{cccc}
							T_0  & &&  \\
							-I_m & T_0&&\\
							&\ddots&\ddots&\\
							&&-I_m& T_0
						\end{array}
						\right],\quad T_0=(1+\alpha)I_m+\tau M_{m}^{-\frac{1}{2}} K_{m} M_{m}^{-\frac{1}{2}},
					\end{equation*}
					and that $(\widetilde{\mathcal{T}}+\alpha I_{mn})^{-1}$ has the following expression
					\begin{align}
						&  (\widetilde{\mathcal{T}}+\alpha I_{mn})^{-1}=\left[
						\begin{array}
							[c]{cccc}
							T_0^{-1}&       &  &\\
							T_0^{-2}& T_0^{-1}&  &\\
							\vdots&\ddots &\ddots&\\
							T_0^{-n}& \ldots & T_0^{-2} & T_0^{-1}
						\end{array}
						\right].\label{tildeTplusalphinvexpress}
					\end{align}
					Denote \correc{with}
					\begin{equation*} \label{eqn:Z_eps}
						\mathcal{Z}_{\epsilon} := \epsilon^{-1} (I_{m} - \epsilon T_0^{-n}).
					\end{equation*}
					\begin{lemma}\label{lem:C_eps_invertible}
						For $\epsilon \in (0,1]$, both $\widetilde{\mathcal{C}}_{\epsilon}+ \alpha {I}_{mn}$ and $\mathcal{Z}_{\epsilon}$ are invertible. Also, we have
						\begin{eqnarray*}
							&&(\widetilde{\mathcal{C}}_{\epsilon}+ \alpha {I}_{mn})^{-1}(\widetilde{\mathcal{T}}+ \alpha {I}_{mn}) =I_{mn} + (\widetilde{\mathcal{T}}+ \alpha {I}_{mn})^{-1}E_{1}\mathcal{Z}_{\epsilon}^{-1} E_{n}^{\top},
						\end{eqnarray*}
                        \correc{where  $E_{i}=e_{i} \otimes I_{m}$ with $e_{i}$ denoting the $i$-th column of $I_{n}$ and $\otimes$ denoting the Kronecker product.}
					\end{lemma}
					\begin{proof}
						It is clear that
						\begin{eqnarray*}
							\sigma(T_0)&=& \sigma \left(\left(\alpha + 1\right)I_{m}+\tau M_{m}^{-\frac{1}{2}} K_{m} M_{m}^{-\frac{1}{2}}\right)\\
							&\subset& (1,+\infty).
						\end{eqnarray*}
						
						Hence, $\sigma( T_0^{-1} ) \subset (0,1)$. Using $\epsilon \in (0,1]$, we have $\sigma(I_{m} - \epsilon T_0^{-n}) \subset (0,1)$, implying $\mathcal{Z}_{\epsilon}$ is invertible.
						
						Then, by \eqref{tildeTplusalphinvexpress}, we see that
						\begin{equation*}
							\epsilon[I_m-\epsilon E_n^{\top}(\widetilde{\mathcal{T}}+ \alpha {I}_{mn})^{-1}E_1]^{-1}=\mathcal{Z}_{\epsilon}^{-1}. 
						\end{equation*}
						
						\correc{Note that $\widetilde{\mathcal{C}}_{\epsilon}+ \alpha {I}_{mn} =\widetilde{\mathcal{T}}+ \alpha {I}_{mn}  - \epsilon E_{1}  E_{n}^{\top}$.}
						By the Sherman–Morrison–Woodbury formula, we have
						\begin{align*}
							(\widetilde{\mathcal{C}}_{\epsilon}+ \alpha {I}_{mn} )^{-1}= (\widetilde{\mathcal{T}}+ \alpha {I}_{mn})^{-1}+  (\widetilde{\mathcal{T}}+ \alpha {I}_{mn})^{-1}E_1\mathcal{Z}_{\epsilon}^{-1}E_n^{\top}(\widetilde{\mathcal{T}}+ \alpha {I}_{mn})^{-1}.
						\end{align*}
						Therefore,
						\begin{equation*}\label{eqn:C_eps_inv}
							(\widetilde{\mathcal{C}}_{\epsilon}+ \alpha {I}_{mn})^{-1}(\widetilde{\mathcal{T}}+ \alpha {I}_{mn})=   I_{mn} +(\widetilde{\mathcal{T}}+ \alpha {I}_{mn})^{-1}E_1\mathcal{Z}_{\epsilon}^{-1}E_n^{\top}.
						\end{equation*}

						The proof is complete.
					\end{proof}
					
					
					\begin{remark}\label{rmk:C_eps_T_invertible}
						Similar to Lemma \ref{lem:C_eps_invertible}, one can also show that the matrix $\widetilde{\mathcal{C}}_{\epsilon}^{\top} + \alpha {I}_{mn}$ is also invertible for $\epsilon \in (0,1]$ and 
						\begin{eqnarray*}\label{eqn:C_eps_T_inv}
							&&(\widetilde{\mathcal{C}}_{\epsilon}^{\top} + \alpha {I}_{mn})^{-1}(\widetilde{\mathcal{T}}^{\top}+ \alpha {I}_{mn}) = I_{mn} + (\widetilde{\mathcal{T}}^{\top}+ \alpha {I}_{mn})^{-1} E_{n}\mathcal{Z}_{\epsilon}^{-1} E_{1}^{\top}.
						\end{eqnarray*}
					\end{remark}
					
					Now, $\widetilde{\mathcal{H}}_{\epsilon}={\rm blockdiag}(\widetilde{\mathcal{C}}_{\epsilon}^{\top} + \alpha {I}_{mn},\widetilde{\mathcal{C}}_{\epsilon} + \alpha {I}_{mn})$, which is clearly invertible.
					
					Denote by $\lambda_{\min}({\bf C})$, the minimum eigenvalue of a Hermitian matrix ${\bf C}$.
					\begin{theorem}\label{thm:eigenvalue_P_eps_inv_A}
						The following statements regarding $\widetilde{\mathcal{H}}_{\epsilon}^{-1}\widetilde{\mathcal{H}}$ hold.
						\begin{enumerate}[(i)]
							\item rank$(\widetilde{\mathcal{H}}_{\epsilon}^{-1}\widetilde{\mathcal{H}}-I_{2mn} )$=$2m$. Furthermore, $\widetilde{\mathcal{H}}_{\epsilon}^{-1}\widetilde{\mathcal{H}}$ has exactly $2(n-1)m$ many eigenvalues equal to $1$.
							\item Given any constant $\eta \in (0,1)$, for $\epsilon \in (0,\eta]$, we have
							\begin{equation*}
								\max_{\lambda\in\sigma(\widetilde{\mathcal{H}}_{\epsilon}^{-1}\widetilde{\mathcal{H}})} | \lambda -1 | \leq \frac{\epsilon}{1-\eta}.
							\end{equation*}
						\end{enumerate}
					\end{theorem}
					\begin{proof}
						Using Lemma \ref{lem:C_eps_invertible} and Remark \ref{rmk:C_eps_T_invertible}, we have
						\begin{align*}
							\widetilde{\mathcal{H}}_{\epsilon}^{-1}\widetilde{\mathcal{H}}-I_{2mn}&=\begin{bmatrix}
								(\widetilde{\mathcal{C}}_{\epsilon}^{\top} + \alpha {I}_{mn})^{-1}(\widetilde{\mathcal{T}}^{\top}+ \alpha {I}_{mn})& \\ &(\widetilde{\mathcal{C}}_{\epsilon}+ \alpha {I}_{mn})^{-1}(\widetilde{\mathcal{T}}+ \alpha {I}_{mn}) 
							\end{bmatrix}-I_{2mn}\\
							&=  \begin{bmatrix}
								(\widetilde{\mathcal{T}}^{\top}+ \alpha {I}_{mn})^{-1} E_{n}\mathcal{Z}_{\epsilon}^{-1} E_{1}^{\top}& \\ &(\widetilde{\mathcal{T}}+ \alpha {I}_{mn})^{-1}E_1\mathcal{Z}_{\epsilon}^{-1}E_n^{\top} 
							\end{bmatrix}.
						\end{align*}
						Since $\textrm{rank}\left((\widetilde{\mathcal{T}}^{\top}+ \alpha {I}_{mn})^{-1} E_{n}\mathcal{Z}_{\epsilon}^{-1} E_{1}^{\top}\right)=m$ and $\textrm{rank}\left((\widetilde{\mathcal{T}}+ \alpha {I}_{mn})^{-1}E_1\mathcal{Z}_{\epsilon}^{-1}E_n^{\top}\right)=m$, we know that $\textrm{rank}\left(\widetilde{\mathcal{H}}_{\epsilon}^{-1}\widetilde{\mathcal{H}}-I_{2mn}\right)=2m$. Moreover, it is straightforward to check that
						\begin{align*}
							(\widetilde{\mathcal{T}}^{\top}+ \alpha {I}_{mn})^{-1} E_{n}\mathcal{Z}_{\epsilon}^{-1} E_{1}^{\top}&=   \correc{\begin{bmatrix}
								T_o^{-n}\mathcal{Z}_{\epsilon}^{-1} & 0 & \cdots  & 0  \\
								T_o^{-n+1}\mathcal{Z}_{\epsilon}^{-1} & \vdots &  & \vdots  \\
								\vdots & \vdots & & \vdots  \\
								T_o^{-1}\mathcal{Z}_{\epsilon}^{-1} & 0 & \cdots & 0
							\end{bmatrix}},\label{resmat1}
						\end{align*}
						\begin{align}
							(\widetilde{\mathcal{T}}+ \alpha {I}_{mn})^{-1}E_1\mathcal{Z}_{\epsilon}^{-1}E_n^{\top}&=    \correc{\begin{bmatrix}
								0 & \cdots & 0 & T_0^{-1}\mathcal{Z}_{\epsilon}^{-1} \\
								\vdots &  & \vdots  &T_0^{-2}\mathcal{Z}_{\epsilon}^{-1} \\
								\vdots  &  & \vdots  & \vdots \\
								0 & \cdots & 0 & T_0^{-n}\mathcal{Z}_{\epsilon}^{-1}
							\end{bmatrix}}.\label{resmat2}
						\end{align}
						
						Then, we see that $\widetilde{\mathcal{H}}_{\epsilon}^{-1}\widetilde{\mathcal{H}}-I_{2mn}$ has exactly $2(n-1)m$ many zero eigenvalues and thus $\widetilde{\mathcal{H}}_{\epsilon}^{-1}\widetilde{\mathcal{H}}$ has exactly $2(n-1)m$ many eigenvalues equal to $1$.

						From the discussion above, we know that 
						\begin{eqnarray*}
							&&\sigma(\widetilde{\mathcal{H}}_{\epsilon}^{-1}\widetilde{\mathcal{H}}) \\
							&=& \{1\} \cup \sigma(I_{mn}+(\widetilde{\mathcal{T}}+ \alpha {I}_{mn})^{-1}E_1\mathcal{Z}_{\epsilon}^{-1}E_n^{\top})\cup\sigma(I_{mn}+ (\widetilde{\mathcal{T}}^{\top}+ \alpha {I}_{mn})^{-1} E_{n}\mathcal{Z}_{\epsilon}^{-1} E_{1}^{\top})\\
							&=&\{1\} \cup \sigma(I_m+T_0^{-n}\mathcal{Z}_{\epsilon}^{-1}).
						\end{eqnarray*}
						
						Then,
						\begin{align*}
							\max_{\lambda\in\sigma(\widetilde{\mathcal{H}}_{\epsilon}^{-1}\widetilde{\mathcal{H}})} | \lambda -1 | &= \max_{\lambda\in\sigma(T_0^{-n}\mathcal{Z}_{\epsilon}^{-1})} | \lambda  |.
						\end{align*}

						Recall that $T_0^{-n}\mathcal{Z}_{\epsilon}^{-1}=   \epsilon T_0^{-n}(I_m-\epsilon T_0^{-n})^{-1}$
						and $T_0=(1+\alpha)I_m+\tau M_{m}^{-\frac{1}{2}} K_{m} M_{m}^{-\frac{1}{2}}$ is a Hermitian positive definite matrix with $\lambda_{\min}(T_0)>1$. Therefore,
						\begin{align*}
							\max_{\lambda\in\sigma(\widetilde{\mathcal{H}}_{\epsilon}^{-1}\widetilde{\mathcal{H}})} | \lambda -1 | &= \max_{\lambda\in\sigma(T_0^{-n}\mathcal{Z}_{\epsilon}^{-1} )} | \lambda|\\
							&=\max_{\lambda\in\sigma(T_0^{-1})} \left|  \frac{\epsilon\lambda^n}{1-\epsilon\lambda^n} \right| \\
							&=\max_{\lambda\in\sigma(T_0^{-1})}  \frac{\epsilon\lambda^n}{1-\epsilon\lambda^n}\leq \sup\limits_{\lambda\in(0,1)} \frac{\epsilon\lambda^n}{1-\epsilon\lambda^n}\leq \frac{\epsilon}{1-\epsilon}\leq \frac{\epsilon}{1-\eta}.
						\end{align*}
						The proof is complete.
					\end{proof}

					\begin{lemma}\label{lem:2_norm_C_inv_T_minus_I}
						Given any constant $\eta \in (0,1)$, for $\epsilon \in (0,\eta]$, we have
						\begin{equation*}
							\| (\widetilde{\mathcal{C}}_{\epsilon}+ \alpha {I}_{mn})^{-1}(\widetilde{\mathcal{T}}+ \alpha {I}_{mn}) - I_{mn} \|_{2} \leq \frac{ \epsilon \sqrt{n}}{1-\eta}.
						\end{equation*}
					\end{lemma}
					\begin{proof}
						From the proof of \eqref{resmat2}, we see that
						\begin{align*}
							(\widetilde{\mathcal{C}}_{\epsilon}+ \alpha {I}_{mn})^{-1}(\widetilde{\mathcal{T}}+ \alpha {I}_{mn}) - I_{mn}&= (\widetilde{\mathcal{T}}+ \alpha {I}_{mn})^{-1}E_1\mathcal{Z}_{\epsilon}^{-1}E_n^{\top}\\
							&=    \correc{\begin{bmatrix}
								0 & \cdots & 0 & T_0^{-1}\mathcal{Z}_{\epsilon}^{-1} \\
								\vdots & &\vdots  &T_0^{-2}\mathcal{Z}_{\epsilon}^{-1} \\
								\vdots&  & \vdots & \vdots \\
								0 & \cdots & 0 & T_0^{-n}\mathcal{Z}_{\epsilon}^{-1}
							\end{bmatrix}}.
						\end{align*}
						Recall that  $T_0=(1+\alpha)I_m+\tau M_{m}^{-\frac{1}{2}} K_{m} M_{m}^{-\frac{1}{2}}$ is Hermitian positive definite.
						Let $T_0^{-1}=\mathcal{Q}_{m} \Lambda \mathcal{Q}_{m}^{\top}$ be the orthogonal diagonalization of $T_0^{-1}$ with $\mathcal{Q}_{m} \in \mathbb{R}^{m\times m}$ being an orthogonal matrix and $\Lambda \in \mathbb{R}^{m\times m}$ being a diagonal matrix with positive diagonal entries. Recall that $\mathcal{Z}_{\epsilon}^{-1}=\epsilon(I_m-\epsilon T_0^{-n})^{-1}$. Then, one can further decompose $(\widetilde{\mathcal{C}}_{\epsilon}+ \alpha {I}_{mn})^{-1}(\widetilde{\mathcal{T}}+ \alpha {I}_{mn}) - I_{mn}$ as follows
						\begin{eqnarray*}
							&&(\widetilde{\mathcal{C}}_{\epsilon}+ \alpha {I}_{mn})^{-1}(\widetilde{\mathcal{T}}+ \alpha {I}_{mn}) - I_{mn}\\
							&=&\epsilon(I_n\otimes \mathcal{Q}_m) \correc{\begin{bmatrix}
								0 & \cdots & 0 & \Lambda^{1} (I_{m} - \epsilon \Lambda^{n})^{-1} \\
								\vdots &  & \vdots & \Lambda^{2} (I_{m} - \epsilon \Lambda^{n})^{-1} \\
								\vdots &  & \vdots & \vdots \\
								0 & \cdots & 0 & \Lambda^{n} (I_{m} - \epsilon \Lambda^{n})^{-1} 
							\end{bmatrix}}(I_n\otimes \mathcal{Q}_m^{\top}).
						\end{eqnarray*}
						Then, by rewriting $\Lambda = {\rm diag}(\lambda_{i})_{i=1}^{m}$, we have
						\begin{eqnarray*}
							&&\|(\widetilde{\mathcal{C}}_{\epsilon}+ \alpha {I}_{mn})^{-1}(\widetilde{\mathcal{T}}+ \alpha {I}_{mn}) - I_{mn} \|_{2} \\
							&\leq& \epsilon \correc{\cdot} \sqrt{\left\| \sum_{j=1}^{n} (\Lambda^{j}(I_{m} - \epsilon \Lambda^{n})^{-1})^{2} \right\|_{2}}\\
							&=& \epsilon \correc{\cdot} \sqrt{\max_{1\leq k \leq m} \sum_{j=1}^{n} \left(\frac{\lambda_{k}^{j}}{1-\epsilon \lambda_{k}^{n}} \right)^{2}}\leq \epsilon \correc{\cdot}  \sqrt{\sup_{\lambda\in(0,1)} \sum_{j=1}^{n} \left(\frac{\lambda^{j}}{1-\epsilon \lambda^{n}} \right)^{2}} ,
						\end{eqnarray*}
						where the last inequality comes from the fact that $\sigma(T_0^{-1})\subset(0,1)$. 
						Since the functions $g_{j}(x) = \frac{x^{j}}{1-\epsilon x^{n}}$ are monotonically increasing on $x\in [0,1]$ for each $j=1,2,\dots,n$,
						we have 
						\begin{eqnarray*}
							\|(\widetilde{\mathcal{C}}_{\epsilon}+ \alpha {I}_{mn})^{-1}(\widetilde{\mathcal{T}}+ \alpha {I}_{mn}) - I_{mn} \|_{2}
							&\leq& \epsilon \correc{\cdot}\sqrt{\sum_{j=1}^{n}\frac{1}{(1-\epsilon)^{2}}}\\
							&=& \frac{\epsilon \sqrt{n}}{1-\epsilon}\\
							&\leq& \frac{ \epsilon \sqrt{n}}{1-\eta}.
						\end{eqnarray*}
						The proof is complete.
					\end{proof}
					
					\begin{remark}\label{rmk:2_norm_C_inv_T_minus_I_transpose}
						\correc{Following  on the proof presented in Lemma~\ref{lem:2_norm_C_inv_T_minus_I}, we can derive the same bound for $\|  (\widetilde{\mathcal{C}}_{\epsilon}^{\top} + \alpha {I}_{mn})^{-1}(\widetilde{\mathcal{T}}^{\top}+ \alpha {I}_{mn}) - I_{mn} \|_{2}$. Specifically, given} any constant $\eta \in (0,1)$, for any $\epsilon \in (0,\eta]$, we have
						\begin{equation*}
							\|  (\widetilde{\mathcal{C}}_{\epsilon}^{\top} + \alpha {I}_{mn})^{-1}(\widetilde{\mathcal{T}}^{\top}+ \alpha {I}_{mn}) - I_{mn} \|_{2} \leq \frac{ \epsilon \sqrt{n}}{1-\eta}.
						\end{equation*}
						
					\end{remark}

					\begin{lemma}\label{lem:2_norm_P_eps_inv_A}
						Given any $\eta \in (0,1)$, choose $\epsilon \in (0,\eta]$. Then,
						\begin{equation*}
							\| \widetilde{\mathcal{H}}_{\epsilon}^{-1}\widetilde{\mathcal{B}}\|_{2} \leq \sqrt{2} \left(1 + \frac{\epsilon \sqrt{n}}{1-\eta}\right).
						\end{equation*}
					\end{lemma}
					\begin{proof}
						\begin{eqnarray*}
							\|  \widetilde{\mathcal{H}}_{\epsilon}^{-1}\widetilde{\mathcal{B}} \|_{2}
							&=&  
							\| \widetilde{\mathcal{H}}_{\epsilon}^{-1}\widetilde{\mathcal{H}} \widetilde{\mathcal{H}} ^{-1} \widetilde{\mathcal{B}} \|_{2} \\
							&\leq&
							\| \widetilde{\mathcal{H}}_{\epsilon}^{-1}\widetilde{\mathcal{H}} \|_{2} \| \widetilde{\mathcal{H}} ^{-1} \widetilde{\mathcal{B}} \|_{2} \\
							&\leq&
							(\| I_{2mn} \|_{2} + \| \widetilde{\mathcal{H}}_{\epsilon}^{-1}\widetilde{\mathcal{H}} - I_{2mn} \|_{2}) \|\widetilde{\mathcal{H}} ^{-1} \widetilde{\mathcal{B}} \|_{2}.
						\end{eqnarray*}
						Since 
						\begin{eqnarray*}
							&&\widetilde{\mathcal{H}}_{\epsilon}^{-1}\widetilde{\mathcal{H}} - I_{2mn}\\
							&=&  \begin{bmatrix}
								(\widetilde{\mathcal{C}}_{\epsilon}^{\top} + \alpha {I}_{mn})^{-1}(\widetilde{\mathcal{T}}^{\top}+ \alpha {I}_{mn}) - I_{mn} & \\
								& (\widetilde{\mathcal{C}}_{\epsilon}+ \alpha {I}_{mn})^{-1}(\widetilde{\mathcal{T}}+ \alpha {I}_{mn}) - I_{mn}
							\end{bmatrix},
						\end{eqnarray*}
						Lemma \ref{lem:2_norm_C_inv_T_minus_I} and Remark \ref{rmk:2_norm_C_inv_T_minus_I_transpose} immediately induce that
						\begin{equation*}
							\| \widetilde{\mathcal{H}}_{\epsilon}^{-1}\widetilde{\mathcal{H}} - I_{2mn} \|_{2} \leq \frac{\epsilon \sqrt{n}}{1-\eta}.
						\end{equation*}
						Therefore,
						\begin{eqnarray*}
							\| \widetilde{\mathcal{H}}_{\epsilon}^{-1}\widetilde{\mathcal{B}} \|_{2}
							&\leq&
							\left(1 + \frac{ \epsilon \sqrt{n}}{1-\eta}\right) \| \widetilde{\mathcal{H}} ^{-1} \widetilde{\mathcal{B}}\|_{2}.
						\end{eqnarray*}
						From Theorem \ref{thm:spectum_of_H_inv_B}, we know that $\sigma(\widetilde{\mathcal{H}} ^{-1} \widetilde{\mathcal{B}}) \subseteq\{1+{\bf i}x|x\in[-1,1]\}$ and that $\widetilde{\mathcal{H}} ^{-1} \widetilde{\mathcal{B}}$ is unitarily diagonalizable, which implies $\|\widetilde{\mathcal{H}} ^{-1} \widetilde{\mathcal{B}}\|_{2} \leq \sqrt{2}$. Therefore,
						\begin{eqnarray*}
							\| \widetilde{\mathcal{H}}_{\epsilon}^{-1}\widetilde{\mathcal{B}}\|_{2}
							&\leq&
							\sqrt{2} \left(1 + \frac{ \epsilon \sqrt{n}}{1-\eta}\right).
						\end{eqnarray*}
						The proof is complete.
					\end{proof}
					
					For any matrix $\mathcal{K} \in \mathbb{R}^{m \times m}$, denote
					\begin{eqnarray*}
						\mathbb{H}(\mathcal{K}) := \frac{\mathcal{K}+\mathcal{K}^{\top}}{2}.
					\end{eqnarray*}
					
					
					\begin{lemma}\label{lem:2_norm_sym_P_eps_inv_A_minus_I}
						Given $\eta \in (0,1)$, for any $\epsilon \in (0,\eta]$, we have 
						\begin{equation*}
							\|\mathbb{H}(\widetilde{\mathcal{H}}_{\epsilon}^{-1}\widetilde{\mathcal{B}}-I_{2mn}) \|_{2} \leq \frac{2  \epsilon\sqrt{n}}{1-\eta}.
						\end{equation*}
						
					\end{lemma}
					\begin{proof}
						Denote \correc{with}
						\begin{equation*}
							\mathcal{R}_1=  (\widetilde{\mathcal{C}}_{\epsilon}^{\top} + \alpha {I}_{mn})^{-1}(\widetilde{\mathcal{T}}^{\top}+ \alpha {I}_{mn}) - I_{mn},\quad \mathcal{R}_2= (\widetilde{\mathcal{C}}_{\epsilon}+ \alpha {I}_{mn})^{-1}(\widetilde{\mathcal{T}}+ \alpha {I}_{mn}) - I_{mn}.
						\end{equation*}
						Recall that $\widetilde{\mathcal{E}}=(\widetilde{\mathcal{T}}+ \alpha I_{mn})^{-1} (-\widetilde{\mathcal{T}} + \alpha I_{mn})$ defined in proof of Theorem \ref{thm:spectum_of_H_inv_B}. 
						
						By the proof of Theorem \ref{thm:spectum_of_H_inv_B} and the proof of Lemma \ref{lem:2_norm_C_inv_T_minus_I}, we know that
						\begin{align*}
							\widetilde{\mathcal{H}}_{\epsilon}^{-1}\widetilde{\mathcal{B}}&=    \widetilde{\mathcal{H}}_{\epsilon}^{-1}\widetilde{\mathcal{H}}\widetilde{\mathcal{H}}^{-1}\widetilde{\mathcal{B}}\\
							&=(I_{2mn}+\widetilde{\mathcal{H}}_{\epsilon}^{-1}\widetilde{\mathcal{H}}-I_{2mn})\begin{bmatrix}
								I_{mn}&-\widetilde{\mathcal{E}}^{\top}\\\widetilde{\mathcal{E}}&I_{mn}
							\end{bmatrix}\\
							&=\begin{bmatrix}
								I_{mn}+\mathcal{R}_1& \\ &I_{mn}+\mathcal{R}_2
							\end{bmatrix}\begin{bmatrix}
								I_{mn}&-\widetilde{\mathcal{E}}^{\top}\\\widetilde{\mathcal{E}}&I_{mn}
							\end{bmatrix}\\
							&=\begin{bmatrix}
								I_{mn}+\mathcal{R}_1& -(I_{mn}+\mathcal{R}_1)\widetilde{\mathcal{E}}^{\top} \\ 
								(I_{mn}+\mathcal{R}_2)\widetilde{\mathcal{E}}&I_{mn}+\mathcal{R}_2 
							\end{bmatrix}.
						\end{align*}

						By simple calculations, we have
						\begin{eqnarray*}
							&&\mathbb{H}(\widetilde{\mathcal{H}}_{\epsilon}^{-1}\widetilde{\mathcal{B}}-I_{2mn})\\
							&=&
							\frac{1}{2}
							\begin{bmatrix}
								\mathcal{R}_1+\correc{\mathcal{R}_1^{\top}} & (\mathcal{R}_2 \widetilde{\mathcal{E}} - \widetilde{\mathcal{E}} \mathcal{R}_1^{\top})^{\top}\\
								\mathcal{R}_2\widetilde{\mathcal{E}}- \widetilde{\mathcal{E}} \mathcal{R}_1^{\top} & \mathcal{R}_2+\mathcal{R}_2^{\top} 
							\end{bmatrix}
							\\
							&=&
							\frac{1}{2}
							\begin{bmatrix}
								\mathcal{R}_1+\correc{\mathcal{R}_1^{\top}} & \\
								& \mathcal{R}_2+\mathcal{R}_2^{\top} 
							\end{bmatrix}
							+ \frac{1}{2}
							\begin{bmatrix}
								&(\mathcal{R}_2\widetilde{\mathcal{E}} - \widetilde{\mathcal{E}}\mathcal{R}_1^{\top})^{\top}\\
								\mathcal{R}_2 \widetilde{\mathcal{E}} - \widetilde{\mathcal{E}} \mathcal{R}_1^{\top}  & 
							\end{bmatrix}.
						\end{eqnarray*}
						Using the result of Lemma \ref{lem:2_norm_C_inv_T_minus_I}, we have
						\begin{eqnarray*} 
							&&\frac{1}{2} \left\| 
							\begin{bmatrix}
								\mathcal{R}_{1}+\mathcal{R}_{1}^{\top} & \\
								& \mathcal{R}_{2}+\mathcal{R}_{2}^{\top}
							\end{bmatrix} \right\|_{2}\\
							&\correc{=}& \frac{1}{2} \max \{\| \mathcal{R}_{1}+\mathcal{R}_{1}^{\top} \|_{2}, \| \mathcal{R}_{2}+\mathcal{R}_{2}^{\top} \|_{2} \}\\
							&\leq& \frac{1}{2} \max \{\| \mathcal{R}_{1} \|_{2} + \| \mathcal{R}_{1}^{\top} \|_{2}, \| \mathcal{R}_{2} \|_{2} + \| \mathcal{R}_{2}^{\top} \|_{2} \}\\
							&=& \max\{\| \mathcal{R}_{1} \|_{2}, \| \mathcal{R}_{2} \|_{2} \}\\
							&\leq& \frac{\epsilon \sqrt{n}}{1-\eta}.
						\end{eqnarray*}
						Note that $\| \widetilde{\mathcal{E}} \|_{2} < 1$ from the result in \cite{bai2014}. Then,
						\correc{\begin{eqnarray*}
							&&\frac{1}{2} \left\|
							\begin{bmatrix}
								& (\mathcal{R}_{2} \widetilde{\mathcal{E}} - \widetilde{\mathcal{E}} \mathcal{R}_{1}^{\top})^{\top}\\
								\mathcal{R}_{2} \widetilde{\mathcal{E}} - \widetilde{\mathcal{E}} \mathcal{R}_{1}^{\top} & 
							\end{bmatrix}
							\right\|_{2}\\
							&=&
							\frac{1}{2} \max\{ \| (\mathcal{R}_{2} \widetilde{\mathcal{E}} - \widetilde{\mathcal{E}} \mathcal{R}_{1}^{\top})^{\top} \|_{2}, \| \mathcal{R}_{2} \widetilde{\mathcal{E}} - \widetilde{\mathcal{E}} \mathcal{R}_{1}^{\top} \|_{2} \}\\
							&=& \frac{1}{2} \| \mathcal{R}_{2} \widetilde{\mathcal{E}} - \widetilde{\mathcal{E}} \mathcal{R}_{1}^{\top} \|_{2}\\
							&\leq& \frac{1}{2} (\| \mathcal{R}_{2} \|_{2} \| \widetilde{\mathcal{E}} \|_{2} + \| \widetilde{\mathcal{E}} \|_{2} \| \mathcal{R}_{1} \|_{2})\\
							&<& \frac{1}{2} \left( \frac{ \epsilon \sqrt{n}}{1-\eta} + \frac{\epsilon \sqrt{n}}{1-\eta} \right)\\
							&=& \frac{ \epsilon \sqrt{n}}{1-\eta}.
						\end{eqnarray*}}
						Therefore, it is clear that
						\begin{equation*}
							\|\mathbb{H}(\widetilde{\mathcal{H}}_{\epsilon}^{-1}\widetilde{\mathcal{B}}-I_{2mn}) \|_{2} \leq \frac{2 \epsilon\sqrt{n}}{1-\eta}.
						\end{equation*}
						The proof is complete.
					\end{proof}
					
					Denote $\mathcal{O}$ the zero matrix of appropriate dimensions. 
					
					\begin{lemma}\cite{Beckermann_2005}\label{lem:gmres_conv_original}
						Let $\Xi \mathbf{q} = \mathbf{w}$ be a real square linear system with $\mathcal{H}(\Xi) \succ \mathcal{O}$. Then, the residuals of the iterates generated by applying GMRES to solving $\Xi \mathbf{q} = \mathbf{w}$ satisfy
						\begin{equation*}
							\|  \mathbf{r_k} \|_{2} \leq \left( 1 - \frac{\lambda_{\textrm{min}}(\mathbb{H}(\Xi))^{2}}{\|\Xi\|_{2}^{2}} \right)^{k/2}  \|  \mathbf{r_0} \|_{2},
						\end{equation*}
						where $\mathbf{r_k}=\mathbf{w}-\Xi\mathbf{q_k}$ with $\mathbf{q_k}~(k\geq 1)$ being the $k$-th GMRES iteration solution and $\mathbf{q_0}$ being an arbitrary initial guess.
					\end{lemma}
					
					\begin{lemma}\label{lm:gmres_conv}
						Given $\delta \in (0,1)$, for any $\epsilon \in (0, c_{\tau}]$, where $c_{\tau} := \frac{\delta\sqrt{\tau}}{\delta\sqrt{\tau}+2  \sqrt{T}}$ \correc{with $\tau=T/(n+1)$}, the residuals of the iterates generated by applying GMRES to solving the auxiliary linear system \eqref{auxiliarysystem} satisfy
						\begin{equation*}
							\|  \tilde{\mathbf{r}}_k \|_{2} \leq \left( \frac{\sqrt{-\delta^2+8\delta+2}}{2+\delta} \right)^{k}  \|  \tilde{\mathbf{r}}_0 \|_{2},
						\end{equation*}
						where $\tilde{\mathbf{r}}_k =\tilde{\mathbf{b}}-\widetilde{\mathcal{H}}_{\epsilon}^{-1}\widetilde{\mathcal{B}}\tilde{\mathbf{x}}_k$ with $\tilde{\mathbf{x}}_k~(k\geq 1)$ being the $k$-th GMRES iteration solution and $\tilde{\mathbf{x}}_0$ being an arbitrary initial guess.
					\end{lemma}
					\begin{proof}
						By Lemma \ref{lem:2_norm_sym_P_eps_inv_A_minus_I}, we have
						\begin{eqnarray*}
							\| \mathbb{H}(\widetilde{\mathcal{H}}_{\epsilon}^{-1}\widetilde{\mathcal{B}}-I_{2mn}) \|_{2} &\leq& \frac{2   \epsilon\sqrt{n}}{1-c_{\tau}} \\
							&=& \frac{\epsilon\delta}{c_{\tau}}\\
							&\leq& \delta.
						\end{eqnarray*}
						Then,
						\begin{eqnarray}\label{eqn:P_eps_inv_A_pos}
							\mathbb{H}(\widetilde{\mathcal{H}}_{\epsilon}^{-1}\widetilde{\mathcal{B}}) &=& I_{2mn} + \correc{\mathbb{H}}(\widetilde{\mathcal{H}}_{\epsilon}^{-1}\widetilde{\mathcal{B}}-I_{2mn})\\\nonumber
                            &\succeq& \correc{I_{2mn}-\| \mathbb{H}(\widetilde{\mathcal{H}}_{\epsilon}^{-1}\widetilde{\mathcal{B}}-I_{2mn}) \|_{2}I_{2mn}}\\
                            \nonumber
							&=& (1-\delta) I_{2mn} \\\nonumber
							&\succ& \mathcal{O}.
						\end{eqnarray}
						Since $  \mathbb{H}(\widetilde{\mathcal{H}}_{\epsilon}^{-1}\widetilde{\mathcal{B}})\succ \mathcal{O}$, Lemma \ref{lem:gmres_conv_original} is applicable to the preconditioned system \eqref{auxiliarysystem}. From (\ref{eqn:P_eps_inv_A_pos}), it is clear that
						\begin{eqnarray*}
							\lambda_{\textrm{min}}\left( \mathbb{H}\left(\widetilde{\mathcal{H}}_{\epsilon}^{-1}\widetilde{\mathcal{B}}\right)\right)^{2} \geq (1-\delta)^{2}.
						\end{eqnarray*}
						
						By Lemma \ref{lem:2_norm_P_eps_inv_A}, we have
						\begin{eqnarray*}
							\|  \mathbb{H}(\widetilde{\mathcal{H}}_{\epsilon}^{-1}\widetilde{\mathcal{B}})\|_{2} &\leq& \sqrt{2} \left(1 + \frac{  \epsilon \sqrt{n}}{1-c_{\tau}}\right)\\
							&\leq& \sqrt{2} \left( 1+\frac{\delta}{2} \right).
						\end{eqnarray*}
						
						Then, Lemma \ref{lem:gmres_conv_original} gives
						\begin{eqnarray*}
							\|  \tilde{\mathbf{r}}_k\|_{2} 
							&\leq&
							\left( 1 - \frac{(1-\delta)^2}{\left(\sqrt{2} \left( 1+\frac{\delta}{2} \right)\right)^2} \right)^{k/2}  \|  \tilde{\mathbf{r}}_0 \|_{2}\\
							&\leq&
							\left( \frac{\sqrt{-\delta^2+8\delta+2}}{2+\delta} \right)^{k}  \|  \tilde{\mathbf{r}}_k\|_{2}.
						\end{eqnarray*}
						The proof is complete.
					\end{proof}
					
					\begin{theorem}\label{thm:gmres_conv}
						Given $\delta \in (0,1)$, for any $\epsilon \in (0, c_{\tau}]$, where $c_{\tau} := \frac{\delta\sqrt{\tau}}{\delta\sqrt{\tau}+2  \sqrt{T}}$ \correc{with $\tau=T/(n+1)$}, the residuals of the iterates generated by applying GMRES to solving the preconditioned linear system \eqref{eqn:preconditioned_system} satisfy
						\begin{equation*}
							\| \mathbf{r}_k \|_{2} \leq \sqrt{c_0}\left( \frac{\sqrt{-\delta^2+8\delta+2}}{2+\delta} \right)^{k}  \|  \mathbf{r}_0 \|_{2},
						\end{equation*}
						where $c_0$ defined in \eqref{eqn:constant_c0} is a positive constant independent of the matrix size, $\mathbf{r}_k =\mathbf{b}-\mathcal{P}_{\epsilon}^{-1}\mathcal{A}\mathbf{x}_k~(k\geq 0)$ with $\mathbf{x}_k~(k\geq 1)$ being the $k$-th GMRES iteration solution and $\mathbf{x}_0$ being an arbitrary initial guess.
					\end{theorem}
					\begin{proof}
						Let $\tilde{\mathbf{x}}_0=\mathcal{W}^{\frac{1}{2}}\mathcal{G}\mathbf{x}_0$. Then, it holds that $\mathbf{x}_0=\mathcal{G}^{-1}\mathcal{W}^{-\frac{1}{2}}\tilde{\mathbf{x}}_0$. Take $\tilde{\mathbf{x}}_0$ as an initial guess of \correc{the GMRES algorithm} for solving the auxiliary linear system \eqref{auxiliarysystem}. Then, according to Theorem \ref{gmrescvgrelationthm}, we see that
						\begin{equation*}
							\left\|\mathbf{r}_{k}\right\|_{2} \leq \sqrt{2}||\mathcal{W}^{-\frac{1}{2}}||_2||\tilde{\bf r}_k||_2,\quad k=1,2,\dots,
						\end{equation*}
						where $\tilde{\mathbf{r}}_{k}:=\tilde{\bf b} - \widetilde{\mathcal{H}}_{\epsilon}^{-1}\widetilde{\mathcal{B}} \tilde{\mathbf{x}}_{k}$ denotes the residual vector at the $k$-th GMRES iteration for solving \eqref{auxiliarysystem} and $\tilde{\mathbf{x}}_{k}$ denotes the $k$-th iterative solution of \correc{GMRES} for solving \eqref{auxiliarysystem}. According to Lemma \ref{lm:gmres_conv}, we can further estimate $||\tilde{\bf r}_k||_2$ as 
						\begin{equation*}
							\|  \tilde{\mathbf{r}}_k \|_{2} \leq \left( \frac{\sqrt{-\delta^2+8\delta+2}}{2+\delta} \right)^{k}  \|  \tilde{\mathbf{r}}_0 \|_{2},
						\end{equation*}
						where $\tilde{\mathbf{r}}_{0}:=\tilde{\bf b} - \widetilde{\mathcal{H}}_{\epsilon}^{-1}\widetilde{\mathcal{B}} \tilde{\mathbf{x}}_{0}$ denotes the initial residual vector. Therefore,
						\begin{equation*}
							\left\|\mathbf{r}_{k}\right\|_{2} \leq \sqrt{2}||\mathcal{W}^{-\frac{1}{2}}||_2\left( \frac{\sqrt{-\delta^2+8\delta+2}}{2+\delta} \right)^{k}  \|  \tilde{\mathbf{r}}_0 \|_{2},\quad k=1,2,\dots. 
						\end{equation*}
						Note that $\tilde{\mathbf{r}}_0=\mathcal{W}^{\frac{1}{2}}\mathcal{G}\mathbf{r}_0$. Therefore,
						\begin{align*}
							\left\|\mathbf{r}_{k}\right\|_{2}& \leq \sqrt{2}||\mathcal{W}^{-\frac{1}{2}}||_2\left( \frac{\sqrt{-\delta^2+8\delta+2}}{2+\delta} \right)^{k}  \|  \mathcal{W}^{\frac{1}{2}}\mathcal{G}\mathbf{r}_0 \|_{2}\\
							&\leq  \sqrt{2}||\mathcal{W}^{-\frac{1}{2}}||_2\left( \frac{\sqrt{-\delta^2+8\delta+2}}{2+\delta} \right)^{k} \|  \mathcal{W}^{\frac{1}{2}} \|_{2} \|  \mathcal{G}\|_{2} \|  \mathbf{r}_0 \|_{2}\\
							&=\kappa_2(\mathcal{W}^{\frac{1}{2}})\left( \frac{\sqrt{-\delta^2+8\delta+2}}{2+\delta} \right)^{k}\|  \mathbf{r}_0 \|_{2}\\
							&=\sqrt{\kappa_2(M_m)}\left( \frac{\sqrt{-\delta^2+8\delta+2}}{2+\delta} \right)^{k}\|  \mathbf{r}_0 \|_{2}\\
							&\leq \sqrt{c_0}\left( \frac{\sqrt{-\delta^2+8\delta+2}}{2+\delta} \right)^{k}\|  \mathbf{r}_0 \|_{2},\quad k=1,2,\dots .
						\end{align*}
						The proof is complete.
					\end{proof}
					
					\begin{remark}
						As a consequence of Theorem \ref{thm:gmres_conv}, GMRES can achieve a mesh-independent convergence rate when $\epsilon = \mathcal{O}(\tau^{1/2})$. \correc{This is because the factor $\frac{\sqrt{-\delta^2 + 8\delta + 2}}{2 + \delta}$ lies in the interval $(0, 1)$ for this choice of $\epsilon$. As the iteration count $k$ increases, $\|\mathbf{r}_k\|_2$ tends to zero uniformly.}
					\end{remark}

					\subsection{Implementation}\label{sub:implementation}
					First, we discuss the computation of $\mathcal{A}\mathbf{v}$ for any given vector $\mathbf{v}$. The computation of \correc{the} matrix-vector product $\mathcal{A}\mathbf{v}$ can be computed in $\mathcal{O}(mn)$ since $\mathcal{A}$ is a sparse matrix consisting of two simple bi-diagonal block Toeplitz matrices. The required storage is of $\mathcal{O}(mn)$. 
					
					
					At each GMRES iteration, the matrix-vector product $\mathcal{P}_{\epsilon}^{-1}\mathbf{v}$ for a given vector $\mathbf{v}$ needs to be computed. Recalling that $\epsilon$-circulant matrices are diagonalizable by the product of a diagonal matrix and a discrete Fourier matrix $\mathbb{F}_{n} =\frac{1}{\sqrt{n}}[\theta_{n}^{(i-1)(j-1)}]_{i,j=1}^{n}\in \mathbb{C}^{n\times n}$ with $\theta_{n} = \exp{(\frac{2\pi \mathbf{i}}{n})}$, we can represent the matrix $C_{\epsilon,n}$ defined by \eqref{eqn:matrix_C_eps_n} using the diagonalization $C_{\epsilon,n}=D_{\epsilon}^{-1} \mathbb{F}_{n} \Lambda_{\epsilon,n} \mathbb{F}_{n}^{*} D_{\epsilon}$\correc{, where $D_{\epsilon}={\rm diag} (\epsilon^{\frac{i-1}{N}})_{i=1}^{n}$}. Note that $\Lambda_{\epsilon,n}$ is a diagonal matrix.
					
					Hence, we can decompose $\mathcal{P}_{\epsilon}$ from (\ref{eqn:matrix_P_eps}) as follows:
					\begin{eqnarray*}
						&&\mathcal{P}_{\epsilon}\\
						&=& \frac{1}{2}\begin{bmatrix} 
							\mathcal{C}^{\top}_{\epsilon}+ \alpha {I}_n \otimes M_m  & \\
							& \mathcal{C}_{\epsilon} + \alpha {I}_n \otimes M_m 
						\end{bmatrix}
						\begin{bmatrix} 
							I_{mn}  &  I_{mn} \\
							-I_{mn}  &  I_{mn} 
						\end{bmatrix}\\
						&=&
						\mathcal{U}
						\begin{bmatrix} 
							( \Lambda^{*}_{\epsilon,n} + \alpha I_{n} ) \otimes M_{m} + \tau I_{n} \otimes \correc{K_m}  & \\
							& ( \Lambda_{\epsilon,n} + \alpha I_{n}) \otimes M_{m} + \tau I_{n} \otimes \correc{K_m} 
						\end{bmatrix}
						\mathcal{U}^{-1}
						\\
						&&\times
						\frac{1}{2}\begin{bmatrix} 
							I_{mn}  &  I_{mn} \\
							-I_{mn}  &  I_{mn} 
						\end{bmatrix}.
					\end{eqnarray*}
					Note that $\mathcal{U}=\begin{bmatrix} 
						(\mathbb{F}_{n}^{*} D_{\epsilon})^{\top} \otimes I_{m} & \\
						& D_{\epsilon}^{-1} \mathbb{F}_{n} \otimes I_{m}
					\end{bmatrix}$ \correc{and} $\Lambda_{\epsilon,n}= {\rm diag}(\lambda_{i-1}^{(\epsilon)})_{i=1}^{n}$ with $\lambda_{k}^{(\epsilon)}=1-\epsilon^{\frac{1}{n}} \theta_{n}^{-k}$.
					
					Therefore, the computation of $\mathbf{w} = \mathcal{P}_{\epsilon}^{-1}\mathbf{v}$ can be implemented by the following four steps.
					
					\begin{enumerate}[1.]
						\item $\textrm{Compute}~\widehat{\mathbf{v}} = \mathcal{U}^{-1}\mathbf{v}$,
						\item $\textrm{Compute}$
						\begin{eqnarray*}
							\widetilde{\mathbf{v}} =\begin{bmatrix} 
								( \Lambda^{*}_{\epsilon,n} + \alpha I_{n} ) \otimes M_{m} + \tau I_{n} \otimes \correc{K_m}  & \\
								& ( \Lambda_{\epsilon,n} + \alpha I_{n}) \otimes M_{m} + \tau I_{n} \otimes \correc{K_m} 
							\end{bmatrix}^{-1}\widehat{\mathbf{v}},
						\end{eqnarray*}
						\item $\textrm{Compute}~\widetilde{\mathbf{w}} = \mathcal{U}\widetilde{\mathbf{v}}$,
						\item $\textrm{Compute}~\mathbf{w} = \left( \frac{1}{2}
						\begin{bmatrix} 
							I_{mn}  &  I_{mn} \\
							-I_{mn}  &  I_{mn} 
						\end{bmatrix} \right)^{-1} \widetilde{\mathbf{w}} = \begin{bmatrix} 
							I_{mn}  &  -I_{mn} \\
							I_{mn}  &  I_{mn} 
						\end{bmatrix} \widetilde{\mathbf{w}}$.
					\end{enumerate}
					
					Both Steps 1 and 3 can be computed by fast Fourier transforms in $\mathcal{O}(mn\log{n})$. In Step 2, the shifted Laplacian systems can be efficiently solved for instance by using the multigrid method. A detailed description of this highly effective implementation can be found in \cite{HeLiu2022} for example. The matrix-vector multiplication in Step 4 requires only $\mathcal{O}(mn)$ operations.
					
					\section{Numerical Result}\label{sec:numerical}
					In this section, we provide numerical results to show the performance of our proposed preconditioners. \correc{All numerical experiments are carried out using MATLAB R2024b on a MacBook Pro equipped with Apple M4 Pro and 24GB RAM.}.
					
					The CPU time in seconds is measured using MATLAB built-in functions $\bold{tic}$ and $\bold{toc}$. All Steps $1 - 3$ in Section \ref{sub:implementation} are implemented by the functions $\bold{dst}$ and $\bold{fft}$ as discrete sine transform and fast Fourier transform respectively. The GMRES solver used is implemented using the built-in \correc{function} $\bold{gmres}$. We choose a zero initial guess and a stopping tolerance of \correc{$10^{-8}$, based on the reduction in relative residual norms observed for the solver tested.}
					
					In the related tables, we denote by `Iter' the number of iterations for solving a linear system by an iterative solver within the given accuracy. Denote \correc{with} `DoF', the number of unknowns in a linear system. \correc{Also, the notation GMRES-$\mathcal{P}_{\epsilon}/\mathcal{P}$ denotes the GMRES solver applied with our proposed preconditioner $\mathcal{P}{_\epsilon}$ or the original RBD preconditioner $\mathcal{P}$ proposed in \cite{bailu2021}.}

                    The time interval $[0,T]$ and the \correc{spatial domain} are partitioned uniformly with the mesh step size \correc{$\tau=T/(n+1)$} and $h=1/(m+1)$, respectively, where $h$ can be found in the \correc{corresponding tables. Note that a finite difference method is used throughout the numerical tests.}
                    
                    \correc{We define the error measure $e_h$ as
					\begin{align*}\label{eqn:error_measure}
						e_h&: = \left\| \begin{bmatrix} \mathbf{y}\\ \mathbf{p} \end{bmatrix} - \begin{bmatrix} y\\ p \end{bmatrix} \right\|_{L^{\infty}_{\tau}(L^{2}(\Omega))}\\
                        &:=h^{d/2}\max\limits_{1\leq j\leq N} \max\left\{||{\bf y}_m^{(j)}-y(\Omega_h,j)||_2,||{\bf p}_m^{(j-1)}-p(\Omega_h,j-1)||_2\right\},
					\end{align*}
                    where ${\bf y}_m^{(j)}$ (${\bf p}_m^{(j-1)}$, respectively) denotes the numerical solution of $y$ ($p$, respectively) at $j$-th ($j-1$-th, respectively) time step; $y(\Omega_h,j)$ ($p(\Omega_h,j-1)$, respectively) denote the values of  analytical solution of  $y$ ($p$, respectively) on the spatial grid at $j$-th ($j-1$-th, respectively) time step.}
					
					For each of the following examples, we take $\epsilon = \min{ \left\{ \frac{1}{2}, \frac{1}{2}\tau \right\}}$ for GMRES-$\mathcal{P}_{\epsilon}$.
					
					
					\begin{example}\label{ex:example1}
						In this example \cite{Linheatopt2022}, we consider the following two-dimensional problem\correc{, which involves solving (\ref{eqn:Cost_functional_heat}) with} $\Omega=(0,1)^2$, $T = 1$, $a(x_1,x_2)=1$, and
						\begin{align*}
							f(x_1,x_2,t)&=(2\pi^2 -1)e^{-t}\sin{(\pi x_1)}\sin{(\pi x_2)},\\
							g(x_1,x_2,t)&=e^{-t}\sin{(\pi x_1)}\sin{(\pi x_2)},
						\end{align*}
						\correc{the} analytical solution of which is given by
						\begin{eqnarray*}
							y(x_1,x_2,t)&=e^{-t}\sin{(\pi x_1)}\sin{(\pi x_2)},\quad p=0.
						\end{eqnarray*}
						
						We remark that $K_m$ can be fast diagonalizable by the discrete sine transform matrix in this example, so we applied the fast sine transforms to solve the shifted \correc{Laplacian linear systems} in Step $2$ of the four-step procedures in Subsection \ref{sub:implementation}.
						
						\correc{Tables \ref{tab:table1_new_precon_revised_part_1} \& \ref{tab:table1_new_precon_revised_part_2}} present the iteration counts, \correc{CPU times, and discretization errors} of GMRES-$\mathcal{P}_{\epsilon}$ when employing the preconditioner with various $\gamma$ values in the backward Euler method. Our observations are as follows: (i) GMRES-$\mathcal{P}_{\epsilon}$ performs excellently and consistently, maintaining stable iteration counts and CPU times across a range of $\gamma$ values; and (ii) the error decreases as the mesh is refined. 

                        \correc{To further illustrate the effectiveness of our proposed solver compared to an existing one, we present the numerical results of GMRES-$\mathcal{P}{\epsilon}$ and GMRES-$\mathcal{P}$ in Tables \ref{tab:table1_new_precon_revised_compare} and \ref{tab:table1_RBD_precon_revised_compare}. Throughout all tests in these tables, the mesh size is fixed at $h = 2^{-5}$. In virtually all cases, as $\tau$ decreases, the corresponding iteration counts for both methods stabilize. Moreover, the CPU time required by our proposed GMRES-$\mathcal{P}{\epsilon}$ is significantly lower when $\tau$ decreases, demonstrating its superior efficiency compared to GMRES-$\mathcal{P}$.}

                        \correc{Figure~\ref{fig:Example_1} presents the plot of the residual ratio, defined as $\frac{\|\mathbf{r}_{j}\|_2}{\|\mathbf{r}_{j-1}\|_2}$, against the iteration number $j$ required for convergence under various regularization parameters. The results demonstrate that the theoretically optimal upper bound derived in Theorem~\ref{thm:gmres_conv}, given by
                        \[
                        \lim_{\delta \to 0^+} \frac{\sqrt{-\delta^2 + 8\delta + 2}}{2 + \delta} = \frac{\sqrt{2}}{2},
                        \]
                        can indeed be achieved. Moreover, in some instances, the residual decreases even more rapidly than this theoretical estimate suggests.}

                        \begin{table}[]
							\caption{Results of GMRES-$\mathcal{P}_{\epsilon}$ for Example \ref{ex:example1} with $\gamma=10^{-10},10^{-8},10^{-6}$.}
							\label{tab:table1_new_precon_revised_part_1}
							\centering
                            \begin{tabular}{ccccccc}
                            \hline
                            \multirow{2}{*}{$\gamma$}   & \multirow{2}{*}{$h$}      & \multirow{2}{*}{$\tau$} & \multirow{2}{*}{DoF} & \multicolumn{3}{c}{GMRES-$\mathcal{P}_{\epsilon}$} \\ \cline{5-7} 
                                                        &                           &                         &                      & Iter         & CPU           & $e_{h,\tau}$        \\ \hline
                            \multirow{9}{*}{$10^{-10}$} & \multirow{3}{*}{$2^{-6}$} & $2^{-6}$                & 508032               & 6            & 0.23          & 2.60e-05            \\
                                                        &                           & $2^{-7}$                & 1016064              & 6            & 0.40          & 9.81e-06            \\
                                                        &                           & $2^{-8}$                & 2032128              & 6            & 0.72          & 3.95e-06            \\ \cline{2-7} 
                                                        & \multirow{3}{*}{$2^{-7}$} & $2^{-6}$                & 2064512              & 6            & 0.81          & 1.94e-05            \\
                                                        &                           & $2^{-7}$                & 4129024              & 6            & 1.65          & 6.11e-06            \\
                                                        &                           & $2^{-8}$                & 8258048              & 6            & 3.48          & 1.97e-06            \\ \cline{2-7} 
                                                        & \multirow{3}{*}{$2^{-8}$} & $2^{-6}$                & 8323200              & 6            & 3.35          & 1.78e-05            \\
                                                        &                           & $2^{-7}$                & 16646400             & 6            & 7.26          & 5.18e-06            \\
                                                        &                           & $2^{-8}$                & 33292800             & 6            & 16.17         & 1.47e-06            \\ \hline
                            \multirow{9}{*}{$10^{-8}$}  & \multirow{3}{*}{$2^{-6}$} & $2^{-6}$                & 508032               & 8            & 0.28          & 2.60e-05            \\
                                                        &                           & $2^{-7}$                & 1016064              & 8            & 0.53          & 9.81e-06            \\
                                                        &                           & $2^{-8}$                & 2032128              & 10           & 1.28          & 3.95e-06            \\ \cline{2-7} 
                                                        & \multirow{3}{*}{$2^{-7}$} & $2^{-6}$                & 2064512              & 8            & 1.04          & 1.94e-05            \\
                                                        &                           & $2^{-7}$                & 4129024              & 8            & 2.07          & 6.11e-06            \\
                                                        &                           & $2^{-8}$                & 8258048              & 10           & 5.36          & 1.97e-06            \\ \cline{2-7} 
                                                        & \multirow{3}{*}{$2^{-8}$} & $2^{-6}$                & 8323200              & 8            & 4.10          & 1.78e-05            \\
                                                        &                           & $2^{-7}$                & 16646400             & 8            & 8.87          & 5.19e-06            \\
                                                        &                           & $2^{-8}$                & 33292800             & 10           & 25.98         & 1.47e-06            \\ \hline
                            \multirow{9}{*}{$10^{-6}$}  & \multirow{3}{*}{$2^{-6}$} & $2^{-6}$                & 508032               & 12           & 0.49          & 2.61e-05            \\
                                                        &                           & $2^{-7}$                & 1016064              & 14           & 1.06          & 9.99e-06            \\
                                                        &                           & $2^{-8}$                & 2032128              & 17           & 2.45          & 4.21e-06            \\ \cline{2-7} 
                                                        & \multirow{3}{*}{$2^{-7}$} & $2^{-6}$                & 2064512              & 12           & 1.54          & 1.95e-05            \\
                                                        &                           & $2^{-7}$                & 4129024              & 14           & 3.71          & 6.22e-06            \\
                                                        &                           & $2^{-8}$                & 8258048              & 17           & 9.66          & 2.10e-06            \\ \cline{2-7} 
                                                        & \multirow{3}{*}{$2^{-8}$} & $2^{-6}$                & 8323200              & 12           & 6.17          & 1.79e-05            \\
                                                        &                           & $2^{-7}$                & 16646400             & 14           & 15.49         & 5.28e-06            \\
                                                        &                           & $2^{-8}$                & 33292800             & 17           & 51.43         & 1.57e-06            \\ \hline
                            \end{tabular}
						\end{table}

                        \begin{table}[]
							\caption{Results of GMRES-$\mathcal{P}_{\epsilon}$ for Example \ref{ex:example1} with $\gamma=10^{-4},10^{-2},1$.}
							\label{tab:table1_new_precon_revised_part_2}
							\centering
                            \begin{tabular}{ccccccc}
\hline
\multirow{2}{*}{$\gamma$}  & \multirow{2}{*}{$h$}      & \multirow{2}{*}{$\tau$} & \multirow{2}{*}{DoF} & \multicolumn{3}{c}{GMRES-$\mathcal{P}_{\epsilon}$} \\ \cline{5-7} 
                           &                           &                         &                      & Iter  & CPU                & $e_{h,\tau}$          \\ \hline
\multirow{9}{*}{$10^{-4}$} & \multirow{3}{*}{$2^{-6}$} & $2^{-6}$                & 508032               & 15    & 0.65  & 3.59e-05  \\
                           &                           & $2^{-7}$                & 1016064              & 13    & 0.96  & 1.92e-05  \\
                           &                           & $2^{-8}$                & 2032128              & 7     & 0.95  & 1.26e-05  \\ \cline{2-7} 
                           & \multirow{3}{*}{$2^{-7}$} & $2^{-6}$                & 2064512              & 15    & 2.03   & 2.69e-05  \\
                           &                           & $2^{-7}$                & 4129024              & 13    & 3.48   & 1.19e-05  \\
                           &                           & $2^{-8}$                & 8258048              & 7     & 3.93   & 6.25e-06  \\ \cline{2-7} 
                           & \multirow{3}{*}{$2^{-8}$} & $2^{-6}$                & 8323200              & 15    & 8.01  & 2.46e-05  \\
                           &                           & $2^{-7}$                & 16646400             & 13    & 14.46   & 1.01e-05  \\
                           &                           & $2^{-8}$                & 33292800             & 7     & 18.77   & 4.68e-06  \\ \hline
\multirow{9}{*}{$10^{-2}$} & \multirow{3}{*}{$2^{-6}$} & $2^{-6}$                & 508032               & 17    & 0.85  & 2.02e-04  \\
                           &                           & $2^{-7}$                & 1016064              & 18    & 1.48   & 1.36e-03  \\
                           &                           & $2^{-8}$                & 2032128              & 18    & 2.63   & 1.03e-03  \\ \cline{2-7} 
                           & \multirow{3}{*}{$2^{-7}$} & $2^{-6}$                & 2064512              & 17    & 2.41   & 1.51e-04  \\
                           &                           & $2^{-7}$                & 4129024              & 18    & 4.78  & 8.49e-05  \\
                           &                           & $2^{-8}$                & 8258048              & 18    & 9.92   & 5.13e-05  \\ \cline{2-7} 
                           & \multirow{3}{*}{$2^{-8}$} & $2^{-6}$                & 8323200              & 17    & 8.94   & 1.38e-04  \\
                           &                           & $2^{-7}$                & 16646400             & 18    & 19.84   & 7.20e-05  \\
                           &                           & $2^{-8}$                & 33292800             & 18    & 52.71   & 3.83e-05  \\ \hline
\multirow{9}{*}{$1$}       & \multirow{3}{*}{$2^{-6}$} & $2^{-6}$                & 508032               & 10    & 0.38  & 2.49e-03  \\
                           &                           & $2^{-7}$                & 1016064              & 10    & 0.69  & 1.68e-03  \\
                           &                           & $2^{-8}$                & 2032128              & 10    & 1.28   & 1.27e-04  \\ \cline{2-7} 
                           & \multirow{3}{*}{$2^{-7}$} & $2^{-6}$                & 2064512              & 10    & 1.25  & 1.86e-04  \\
                           &                           & $2^{-7}$                & 4129024              & 10    & 2.62   & 1.04e-04  \\
                           &                           & $2^{-8}$                & 8258048              & 10    & 5.44   & 6.31e-05  \\ \cline{2-7} 
                           & \multirow{3}{*}{$2^{-8}$} & $2^{-6}$                & 8323200              & 10    & 5.28   & 1.70e-04  \\
                           &                           & $2^{-7}$                & 16646400             & 10    & 11.21   & 8.86e-05  \\
                           &                           & $2^{-8}$                & 33292800             & 10    & 25.21   & 4.72e-05  \\ \hline
\end{tabular}
						\end{table}

\begin{table}[]
\caption{Results of GMRES-$\mathcal{P}_{\epsilon}$ for Example \ref{ex:example1} with various regularization parameters $\gamma$}
\label{tab:table1_new_precon_revised_compare}
\centering
\begin{tabular}{ccccccc}
\hline
\multirow{2}{*}{$\gamma$} & \multirow{2}{*}{$h$} & \multirow{2}{*}{$\tau$} & \multirow{2}{*}{DoF} & \multicolumn{3}{c}{GMRES-$\mathcal{P}_{\epsilon}$} \\ \cline{5-7}
& & & & Iter & CPU & $e_{h,\tau}$ \\ \hline

\multirow{7}{*}{$10^{-10}$} & \multirow{7}{*}{$2^{-5}$}
& $2^{-5}$  & 61504 & 6 & 0.05 & 1.13e-04 \\
& & $2^{-6}$  & 123008 & 6 & 0.08 & 5.21e-05 \\
& & $2^{-7}$  & 246016 & 6 & 0.15 & 2.46e-05 \\
& & $2^{-8}$  & 492032 & 6 & 0.24 & 1.19e-05 \\
& & $2^{-9}$  & 984064 & 6 & 0.44 & 5.82e-06 \\
& & $2^{-10}$ & 1968128 & 8 & 0.98 & 2.88e-06 \\
& & $2^{-11}$ & 3936256 & 8 & 2.12 & 1.43e-06 \\ \hline

\multirow{7}{*}{$10^{-8}$} & \multirow{7}{*}{$2^{-5}$}
& $2^{-5}$  & 61504 & 6 & 0.06 & 1.13e-04 \\
& & $2^{-6}$  & 123008 & 8 & 0.16 & 5.21e-05 \\
& & $2^{-7}$  & 246016 & 8 & 0.23 & 2.46e-05 \\
& & $2^{-8}$  & 492032 & 10 & 0.47 & 1.19e-05 \\
& & $2^{-9}$  & 984064 & 11 & 0.84 & 5.84e-06 \\
& & $2^{-10}$ & 1968128 & 13 & 1.78 & 2.91e-06 \\
& & $2^{-11}$ & 3936256 & 16 & 4.20 & 1.49e-06 \\ \hline

\multirow{7}{*}{$10^{-6}$} & \multirow{7}{*}{$2^{-5}$}
& $2^{-5}$  & 61504 & 10 & 0.09 & 1.13e-04 \\
& & $2^{-6}$  & 123008 & 12 & 0.22 & 5.24e-05 \\
& & $2^{-7}$  & 246016 & 14 & 0.43 & 2.51e-05 \\
& & $2^{-8}$  & 492032 & 17 & 0.89 & 1.27e-05 \\
& & $2^{-9}$  & 984064 & 18 & 1.43 & 7.12e-06 \\
& & $2^{-10}$ & 1968128 & 17 & 2.27 & 4.75e-06 \\
& & $2^{-11}$ & 3936256 & 15 & 3.92 & 3.75e-06 \\ \hline

\multirow{7}{*}{$10^{-4}$} & \multirow{7}{*}{$2^{-5}$}
& $2^{-5}$  & 61504 & 15 & 0.16 & 1.29e-04 \\
& & $2^{-6}$  & 123008 & 15 & 0.32 & 7.20e-05 \\
& & $2^{-7}$  & 246016 & 13 & 0.42 & 4.81e-05 \\
& & $2^{-8}$  & 492032 & 7 & 0.30 & 3.78e-05 \\
& & $2^{-9}$  & 984064 & 12 & 0.95 & 3.32e-05 \\
& & $2^{-10}$ & 1968128 & 13 & 1.79 & 3.10e-05 \\
& & $2^{-11}$ & 3936256 & 15 & 3.92 & 2.99e-05 \\ \hline

\multirow{7}{*}{$10^{-2}$} & \multirow{7}{*}{$2^{-5}$}
& $2^{-5}$  & 61504 & 16 & 0.20 & 5.30e-04 \\
& & $2^{-6}$  & 123008 & 17 & 0.33 & 4.06e-04 \\
& & $2^{-7}$  & 246016 & 18 & 0.64 & 3.42e-04 \\
& & $2^{-8}$  & 492032 & 18 & 0.95 & 3.10e-04 \\
& & $2^{-9}$  & 984064 & 18 & 1.43 & 2.94e-04 \\
& & $2^{-10}$ & 1968128 & 19 & 2.75 & 2.86e-04 \\
& & $2^{-11}$ & 3936256 & 19 & 5.11 & 2.82e-04 \\ \hline

\multirow{7}{*}{$1$} & \multirow{7}{*}{$2^{-5}$}
& $2^{-5}$  & 61504 & 10 & 0.12 & 6.52e-04 \\
& & $2^{-6}$  & 123008 & 10 & 0.18 & 4.99e-04 \\
& & $2^{-7}$  & 246016 & 10 & 0.27 & 4.21e-04 \\
& & $2^{-8}$  & 492032 & 10 & 0.46 & 3.82e-04 \\
& & $2^{-9}$  & 984064 & 10 & 0.77 & 3.62e-04 \\
& & $2^{-10}$ & 1968128 & 10 & 1.47 & 3.52e-04 \\
& & $2^{-11}$ & 3936256 & 10 & 2.96 & 3.47e-04 \\ \hline

\end{tabular}
\end{table}

\begin{table}[]
\caption{Results of GMRES-$\mathcal{P}$ for Example \ref{ex:example1} with various regularization parameters $\gamma$}
\label{tab:table1_RBD_precon_revised_compare}
\centering
\begin{tabular}{ccccccc}
\hline
\multirow{2}{*}{$\gamma$} & \multirow{2}{*}{$h$} & \multirow{2}{*}{$\tau$} & \multirow{2}{*}{DoF} & \multicolumn{3}{c}{GMRES-$\mathcal{P}$} \\ \cline{5-7}
& & & & Iter & CPU & $e_{h,\tau}$ \\ \hline

\multirow{7}{*}{$10^{-10}$} & \multirow{7}{*}{$2^{-5}$}
& $2^{-5}$  & 61504    & 6  & 0.14 & 1.13e-04 \\
& & $2^{-6}$  & 123008   & 6  & 0.19 & 5.21e-05 \\
& & $2^{-7}$  & 246016   & 6  & 0.37 & 2.46e-05 \\
& & $2^{-8}$  & 492032   & 6  & 0.70 & 1.19e-05 \\
& & $2^{-9}$  & 984064   & 6  & 1.36 & 5.82e-06 \\
& & $2^{-10}$ & 1968128  & 8  & 3.59 & 2.88e-06 \\
& & $2^{-11}$ & 3936256  & 8  & 7.27 & 1.43e-06 \\ \hline

\multirow{7}{*}{$10^{-8}$} & \multirow{7}{*}{$2^{-5}$}
& $2^{-5}$  & 61504    & 6  & 0.13 & 1.13e-04 \\
& & $2^{-6}$  & 123008   & 8  & 0.28 & 5.21e-05 \\
& & $2^{-7}$  & 246016   & 8  & 0.52 & 2.46e-05 \\
& & $2^{-8}$  & 492032   & 10 & 1.18 & 1.19e-05 \\
& & $2^{-9}$  & 984064   & 11 & 2.53 & 5.84e-06 \\
& & $2^{-10}$ & 1968128  & 13 & 6.08 & 2.91e-06 \\
& & $2^{-11}$ & 3936256  & 16 & 13.53 & 1.49e-06 \\ \hline

\multirow{7}{*}{$10^{-6}$} & \multirow{7}{*}{$2^{-5}$}
& $2^{-5}$  & 61504    & 10 & 0.17 & 1.13e-04 \\
& & $2^{-6}$  & 123008   & 12 & 0.41 & 5.24e-05 \\
& & $2^{-7}$  & 246016   & 14 & 0.92 & 2.51e-05 \\
& & $2^{-8}$  & 492032   & 17 & 2.28 & 1.27e-05 \\
& & $2^{-9}$  & 984064   & 18 & 4.45 & 7.12e-06 \\
& & $2^{-10}$ & 1968128  & 17 & 7.41 & 4.75e-06 \\
& & $2^{-11}$ & 3936256  & 15 & 12.96 & 3.75e-06 \\ \hline

\multirow{7}{*}{$10^{-4}$} & \multirow{7}{*}{$2^{-5}$}
& $2^{-5}$  & 61504    & 15 & 0.29 & 1.29e-04 \\
& & $2^{-6}$  & 123008   & 15 & 0.54 & 7.20e-05 \\
& & $2^{-7}$  & 246016   & 13 & 0.91 & 4.81e-05 \\
& & $2^{-8}$  & 492032   & 7  & 0.90 & 3.78e-05 \\
& & $2^{-9}$  & 984064   & 12 & 3.06 & 3.32e-05 \\
& & $2^{-10}$ & 1968128  & 13 & 6.16 & 3.10e-05 \\
& & $2^{-11}$ & 3936256  & 15 & 13.52 & 2.99e-05 \\ \hline

\multirow{7}{*}{$10^{-2}$} & \multirow{7}{*}{$2^{-5}$}
& $2^{-5}$  & 61504    & 16 & 0.31 & 5.30e-04 \\
& & $2^{-6}$  & 123008   & 17 & 0.66 & 4.06e-04 \\
& & $2^{-7}$  & 246016   & 18 & 1.28 & 3.42e-04 \\
& & $2^{-8}$  & 492032   & 18 & 2.53 & 3.10e-04 \\
& & $2^{-9}$  & 984064   & 18 & 4.70 & 2.94e-04 \\
& & $2^{-10}$ & 1968128  & 19 & 9.20 & 2.86e-04 \\
& & $2^{-11}$ & 3936256  & 19 & 17.38 & 2.82e-04 \\ \hline

\multirow{7}{*}{$1$} & \multirow{7}{*}{$2^{-5}$}
& $2^{-5}$  & 61504    & 10 & 0.18 & 6.52e-04 \\
& & $2^{-6}$  & 123008   & 10 & 0.38 & 4.99e-04 \\
& & $2^{-7}$  & 246016   & 10 & 0.66 & 4.21e-04 \\
& & $2^{-8}$  & 492032   & 10 & 1.25 & 3.82e-04 \\
& & $2^{-9}$  & 984064   & 10 & 2.44 & 3.62e-04 \\
& & $2^{-10}$ & 1968128  & 10 & 5.07 & 3.52e-04 \\
& & $2^{-11}$ & 3936256  & 10 & 10.06 & 3.47e-04 \\ \hline

\end{tabular}
\end{table}

\begin{figure}[h!]
    \centering
    \includegraphics[scale=0.75]{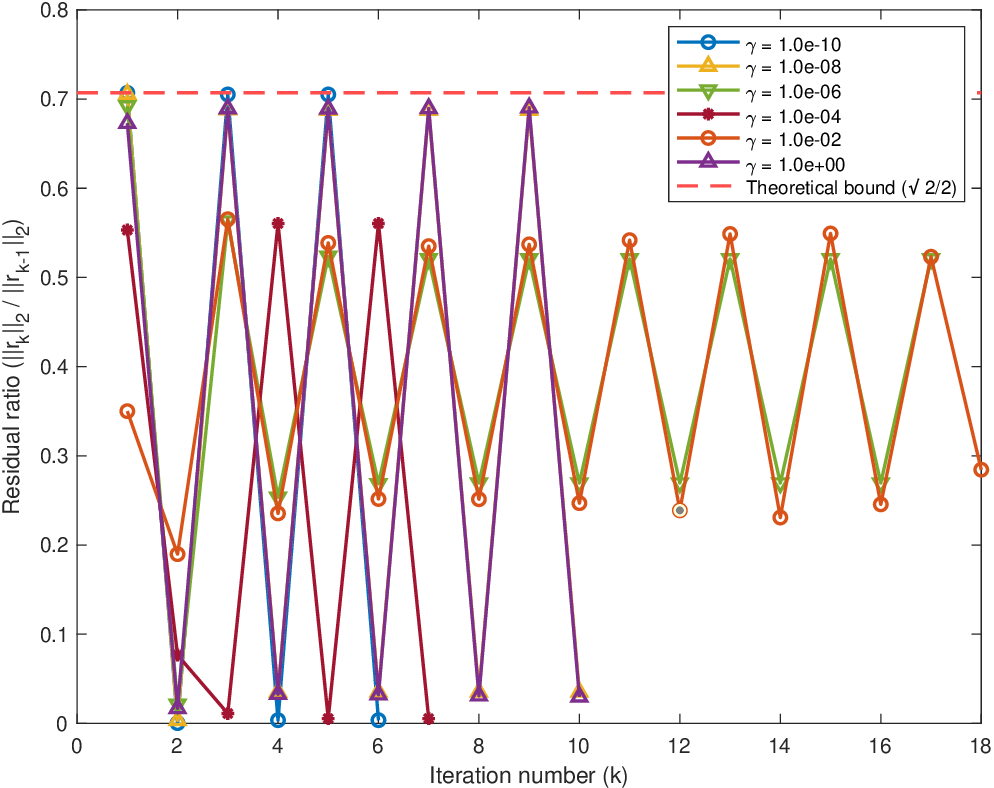}

        \caption{Convergence factor per iteration for Example~\ref{ex:example1} with various regularization parameters $\gamma$.}
        
    \label{fig:Example_1}
\end{figure}

					\end{example}
					
					\begin{example}\label{ex:example2}
						In this example, we consider the following two-dimensional problem\correc{, which involves solving (\ref{eqn:Cost_functional_heat}) with a variable function $a(x_1,x_2)$, $\Omega=(0,1)^2$}, $T=1$, $a(x_1,x_2)=10^{-5} \sin{(\pi x_1 x_2)}$, and
						\begin{multline*}
							f(x_1,x_2,t)= - \sin{(\pi t)}\sin{(\pi x_1)}\sin{(\pi x_2)} + e^{-t} x_1 (1-x_1) [2\times 10^{-5}\sin{(\pi x_1 x_2)}\\
							- x_2 (1-x_2) - 10^{-5} \pi \cos{(\pi x_1 x_2)} x (1-2x_2)]\\
							+e^{-t} x_2 (1-x_2) [2\times 10^{-5} \sin{(\pi x_1 x_2)} - 10^{-5} \pi \cos{(\pi x_1 x_2)} x_2 (1-2x_1)],
						\end{multline*}
						\begin{multline*}
							g(x_1,x_2,t)=-\gamma \pi \cos{(\pi t)}\sin{(\pi x_1)}\sin{(\pi x_2)} + e^{-t} x_1 (1-x_1) x_2 (1-x_2)\\
							- 10^{-5}\gamma \pi^2 \sin{(\pi t)}[-2\sin{(\pi x_1 x_2)}\sin{(\pi x_1)}\sin{(\pi x_2)} \\
							+ \cos{(\pi x_1 x_2)}(x_1\sin{(\pi x_1)}\cos{(\pi x_2)}+x_2\cos{(\pi x_1)}\sin{(\pi x_2)})].
						\end{multline*}
						The analytical solution of which is given by
						\begin{eqnarray*}
							y(x_1,x_2,t)=e^{-t} x_1 (1-x_1) x_2 (1-x_2),\quad p(x_1,x_2,t)=\gamma\sin{(\pi t)}\sin{(\pi x_1)}\sin{(\pi x_2)}.
						\end{eqnarray*}

						Since $K_m$ is not fast diagonalizable in this example, we applied one iteration of the V-cycle geometric multigrid method to solve the shifted Laplacian linear system (as detailed in Subsection \ref{sub:implementation}). The Gauss-Seidel method is employed as the pre-smoother for the multigrid method\correc{, and no post-smoother is used.}
						
						\correc{Tables \ref{tab:table2_new_precon_revised_part_1} \& \ref{tab:table2_new_precon_revised_part_2}} present the iteration counts, \correc{CPU times, and discretization errors} for GMRES-$\mathcal{P}_{\epsilon}$ when the preconditioner is utilized with various values of $\gamma$ in the backward Euler method. The purpose of this example is to evaluate the effectiveness of our solvers when the coefficient function $a(x_1,x_2)$ is non-constant.

                        \correc{Similar to the previous example, we present the numerical results of GMRES-$\mathcal{P}{\epsilon}$ and GMRES-$\mathcal{P}$ in Tables~\ref{tab:table2_new_precon_revised_compare} and~\ref{tab:table2_RBD_precon_revised_compare} for comparison purposes. Throughout all tests reported in these tables, the mesh size is fixed at $h = 2^{-5}$. In nearly all cases, as the time-step size $\tau$ decreases, the iteration counts for both methods tend to stabilize. Furthermore, the CPU time required by our proposed GMRES-$\mathcal{P}_{\epsilon}$ solver is consistently lower for smaller values of $\tau$, showcasing better computational performance than GMRES-$\mathcal{P}$.}
						
						Similar to the last example, the results indicate that (i) GMRES-$\mathcal{P}_{\epsilon}$ achieves stable iteration counts and CPU times across a broad range of $\gamma$ values; (ii) \correc{the error decreases} as expected as the mesh refines.

                        \correc{Similarly to Figure~\ref{fig:Example_1}, Figure~\ref{fig:Example_2} illustrates the residual ratio $\frac{\|\mathbf{r}_{j}\|_2}{\|\mathbf{r}_{j-1}\|_2}$ plotted against the iteration number $j$ for various regularization parameters. The convergence behavior remains consistent with the theoretical upper bound from Theorem~\ref{thm:gmres_conv}, where the optimal upper bound $\frac{\sqrt{2}}{2}$ is attained. In fact, as observed previously, the residual reduction occasionally outpaces this theoretical estimate.
}

\begin{table}[]
\caption{Results of GMRES-$\mathcal{P}_{\epsilon}$ for Example \ref{ex:example2} with $\gamma=10^{-10},10^{-8},10^{-6}$}
\label{tab:table2_new_precon_revised_part_1}
\centering
\begin{tabular}{ccccccc}
\hline
\multirow{2}{*}{$\gamma$} & \multirow{2}{*}{$h$} & \multirow{2}{*}{$\tau$} & \multirow{2}{*}{DoF} & \multicolumn{3}{c}{GMRES-$\mathcal{P}_{\epsilon}$} \\ \cline{5-7}
& & & & Iter & CPU & $e_{h,\tau}$ \\ \hline
\multirow{9}{*}{$10^{-10}$} & \multirow{3}{*}{$2^{-6}$} & $2^{-6}$ & 508032 & 6 & 0.35 & 1.50e-06 \\
& & $2^{-7}$ & 1016064 & 6 & 0.66 & 3.75e-07 \\
& & $2^{-8}$ & 2032128 & 6 & 1.30 & 9.37e-08 \\ \cline{2-7}
& \multirow{3}{*}{$2^{-7}$} & $2^{-6}$ & 2064512 & 6 & 1.43 & 1.50e-06 \\
& & $2^{-7}$ & 4129024 & 6 & 2.58 & 3.75e-07 \\
& & $2^{-8}$ & 8258048 & 6 & 5.01 & 9.37e-08 \\ \cline{2-7}
& \multirow{3}{*}{$2^{-8}$} & $2^{-6}$ & 8323200 & 6 & 5.14 & 1.50e-06 \\
& & $2^{-7}$ & 16646400 & 6 & 10.36 & 3.75e-07 \\
& & $2^{-8}$ & 33292800 & 6 & 21.24 & 9.37e-08 \\ \hline
\multirow{9}{*}{$10^{-8}$} & \multirow{3}{*}{$2^{-6}$} & $2^{-6}$ & 508032 & 8 & 0.39 & 1.50e-06 \\
& & $2^{-7}$ & 1016064 & 8 & 0.77 & 3.75e-07 \\
& & $2^{-8}$ & 2032128 & 10 & 1.95 & 9.37e-08 \\ \cline{2-7}
& \multirow{3}{*}{$2^{-7}$} & $2^{-6}$ & 2064512 & 8 & 1.70 & 1.50e-06 \\
& & $2^{-7}$ & 4129024 & 8 & 3.25 & 3.75e-07 \\
& & $2^{-8}$ & 8258048 & 10 & 7.76 & 9.37e-08 \\ \cline{2-7}
& \multirow{3}{*}{$2^{-8}$} & $2^{-6}$ & 8323200 & 8 & 6.38 & 1.50e-06 \\
& & $2^{-7}$ & 16646400 & 8 & 12.87 & 3.75e-07 \\
& & $2^{-8}$ & 33292800 & 10 & 33.30 & 9.37e-08 \\ \hline
\multirow{9}{*}{$10^{-6}$} & \multirow{3}{*}{$2^{-6}$} & $2^{-6}$ & 508032 & 12 & 0.62 & 1.51e-06 \\
& & $2^{-7}$ & 1016064 & 15 & 1.52 & 3.81e-07 \\
& & $2^{-8}$ & 2032128 & 18 & 3.69 & 9.94e-08 \\ \cline{2-7}
& \multirow{3}{*}{$2^{-7}$} & $2^{-6}$ & 2064512 & 12 & 2.54 & 1.51e-06 \\
& & $2^{-7}$ & 4129024 & 15 & 5.94 & 3.81e-07 \\
& & $2^{-8}$ & 8258048 & 18 & 13.32 & 9.94e-08 \\ \cline{2-7}
& \multirow{3}{*}{$2^{-8}$} & $2^{-6}$ & 8323200 & 12 & 9.27 & 1.51e-06 \\
& & $2^{-7}$ & 16646400 & 15 & 23.49 & 3.81e-07 \\
& & $2^{-8}$ & 33292800 & 18 & 58.45 & 9.94e-08 \\ \hline
\end{tabular}
\end{table}

\begin{table}[]
\caption{Results of GMRES-$\mathcal{P}_{\epsilon}$ for Example \ref{ex:example2} with $\gamma=10^{-4},10^{-2},1$}
\label{tab:table2_new_precon_revised_part_2}
\centering
\begin{tabular}{ccccccc}
\hline
\multirow{2}{*}{$\gamma$}  & \multirow{2}{*}{$h$}      & \multirow{2}{*}{$\tau$} & \multirow{2}{*}{DoF} & \multicolumn{3}{c}{GMRES-$\mathcal{P}_{\epsilon}$} \\ \cline{5-7}
&                           &                           &                          & Iter  & CPU       & $e_{h,\tau}$     \\ \hline
\multirow{9}{*}{$10^{-4}$} & \multirow{3}{*}{$2^{-6}$}  & $2^{-6}$                & 508032   & 20 & 1.12 & 3.84e-06 \\
&                           & $2^{-7}$                & 1016064  & 18 & 1.92 & 1.92e-06 \\
&                           & $2^{-8}$                & 2032128  & 15 & 3.12 & 9.59e-07 \\ \cline{2-7}
& \multirow{3}{*}{$2^{-7}$} & $2^{-6}$                & 2064512  & 20 & 4.46 & 3.84e-06 \\
&                           & $2^{-7}$                & 4129024  & 18 & 7.60 & 1.92e-06 \\
&                           & $2^{-8}$                & 8258048  & 15 & 11.57 & 9.59e-07 \\ \cline{2-7}
& \multirow{3}{*}{$2^{-8}$} & $2^{-6}$                & 8323200  & 20 & 16.27 & 3.84e-06 \\
&                           & $2^{-7}$                & 16646400 & 18 & 29.74 & 1.92e-06 \\
&                           & $2^{-8}$                & 33292800 & 15 & 49.14 & 9.59e-07 \\ \hline
\multirow{9}{*}{$10^{-2}$} & \multirow{3}{*}{$2^{-6}$}  & $2^{-6}$                & 508032   & 12 & 0.61 & 3.50e-04 \\
&                           & $2^{-7}$                & 1016064  & 11  & 1.19 & 1.75e-04 \\
&                           & $2^{-8}$                & 2032128  & 10  & 2.13 & 8.75e-05 \\ \cline{2-7}
& \multirow{3}{*}{$2^{-7}$} & $2^{-6}$                & 2064512  & 12 & 2.52 & 3.50e-04 \\
&                           & $2^{-7}$                & 4129024  & 11 & 4.78 & 1.75e-04 \\
&                           & $2^{-8}$                & 8258048  & 10 & 8.01 & 8.75e-05 \\ \cline{2-7}
& \multirow{3}{*}{$2^{-8}$} & $2^{-6}$                & 8323200  & 12 & 9.57 & 3.50e-04 \\
&                           & $2^{-7}$                & 16646400 & 11 & 18.97 & 1.75e-04 \\
&                           & $2^{-8}$                & 33292800 & 10 & 33.68 & 8.75e-05 \\ \hline
\multirow{9}{*}{$1$}       & \multirow{3}{*}{$2^{-6}$}  & $2^{-6}$                & 508032   & 10 & 0.52 & 1.85e-02 \\
&                           & $2^{-7}$                & 1016064  & 8  & 0.94 & 9.24e-03 \\
&                           & $2^{-8}$                & 2032128  & 8  & 1.87 & 4.61e-03 \\ \cline{2-7}
& \multirow{3}{*}{$2^{-7}$} & $2^{-6}$                & 2064512  & 10 & 2.28 & 1.85e-02 \\
&                           & $2^{-7}$                & 4129024  & 8  & 3.75 & 9.24e-03 \\
&                           & $2^{-8}$                & 8258048  & 8  & 7.23 & 4.61e-03 \\ \cline{2-7}
& \multirow{3}{*}{$2^{-8}$} & $2^{-6}$                & 8323200  & 10 & 8.16 & 1.85e-02 \\
&                           & $2^{-7}$                & 16646400 & 8  & 14.61 & 9.24e-03 \\
&                           & $2^{-8}$                & 33292800 & 8  & 29.12 & 4.61e-03 \\ \hline
\end{tabular}
\end{table}

\begin{table}[]
\caption{Results of GMRES-$\mathcal{P}_{\epsilon}$ for Example \ref{ex:example2} with various regularization parameters $\gamma$}
\label{tab:table2_new_precon_revised_compare}
\centering
\begin{tabular}{ccccccc}
\hline
\multirow{2}{*}{$\gamma$} & \multirow{2}{*}{$h$} & \multirow{2}{*}{$\tau$} & \multirow{2}{*}{DoF} & \multicolumn{3}{c}{GMRES-$\mathcal{P}_{\epsilon}$} \\ \cline{5-7}
& & & & Iter & CPU & $e_{h,\tau}$ \\ \hline
\multirow{7}{*}{$10^{-10}$} & \multirow{7}{*}{$2^{-5}$} & $2^{-5}$ & 61504 & 6 & 0.07 & 6.05e-06 \\
& & $2^{-6}$ & 123008 & 6 & 0.12 & 1.50e-06 \\
& & $2^{-7}$ & 246016 & 6 & 0.24 & 3.75e-07 \\
& & $2^{-8}$ & 492032 & 6 & 0.36 & 9.37e-08 \\
& & $2^{-9}$ & 984064 & 6 & 0.65 & 2.34e-08 \\
& & $2^{-10}$ & 1968128 & 8 & 1.66 & 5.85e-09 \\
& & $2^{-11}$ & 3936256 & 8 & 3.02 & 4.93e-09 \\ \hline
\multirow{7}{*}{$10^{-8}$} & \multirow{7}{*}{$2^{-5}$} & $2^{-5}$ & 61504 & 6 & 0.06 & 6.05e-06 \\
& & $2^{-6}$ & 123008 & 8 & 0.15 & 1.50e-06 \\
& & $2^{-7}$ & 246016 & 8 & 0.29 & 3.75e-07 \\
& & $2^{-8}$ & 492032 & 10 & 0.61 & 9.37e-08 \\
& & $2^{-9}$ & 984064 & 11 & 1.21 & 2.35e-08 \\
& & $2^{-10}$ & 1968128 & 13 & 2.83 & 5.91e-09 \\
& & $2^{-11}$ & 3936256 & 16 & 6.03 & 1.66e-09 \\ \hline
\multirow{7}{*}{$10^{-6}$} & \multirow{7}{*}{$2^{-5}$} & $2^{-5}$ & 61504 & 10 & 0.09 & 6.05e-06 \\
& & $2^{-6}$ & 123008 & 12 & 0.20 & 1.51e-06 \\
& & $2^{-7}$ & 246016 & 15 & 0.50 & 3.81e-07 \\
& & $2^{-8}$ & 492032 & 18 & 1.03 & 9.94e-08 \\
& & $2^{-9}$ & 984064 & 19 & 1.92 & 2.84e-08 \\
& & $2^{-10}$ & 1968128 & 17 & 3.46 & 9.56e-09 \\
& & $2^{-11}$ & 3936256 & 15 & 5.62 & 3.71e-09 \\ \hline
\multirow{7}{*}{$10^{-4}$} & \multirow{7}{*}{$2^{-5}$} & $2^{-5}$ & 61504 & 19 & 0.17 & 7.67e-06 \\
& & $2^{-6}$ & 123008 & 20 & 0.35 & 3.84e-06 \\
& & $2^{-7}$ & 246016 & 18 & 0.63 & 1.92e-06 \\
& & $2^{-8}$ & 492032 & 15 & 0.89 & 9.59e-07 \\
& & $2^{-9}$ & 984064 & 13 & 1.29 & 4.79e-07 \\
& & $2^{-10}$ & 1968128 & 11 & 2.15 & 2.40e-07 \\
& & $2^{-11}$ & 3936256 & 9 & 3.76 & 1.20e-07 \\ \hline
\multirow{7}{*}{$10^{-2}$} & \multirow{7}{*}{$2^{-5}$} & $2^{-5}$ & 61504 & 13 & 0.11 & 7.00e-04 \\
& & $2^{-6}$ & 123008 & 12 & 0.22 & 3.50e-04 \\
& & $2^{-7}$ & 246016 & 11 & 0.39 & 1.75e-04 \\
& & $2^{-8}$ & 492032 & 10 & 0.55 & 8.75e-05 \\
& & $2^{-9}$ & 984064 & 10 & 0.97 & 4.38e-05 \\
& & $2^{-10}$ & 1968128 & 8 & 1.72 & 2.19e-05 \\
& & $2^{-11}$ & 3936256 & 8 & 3.25 & 1.09e-05 \\ \hline
\multirow{7}{*}{$1$} & \multirow{7}{*}{$2^{-5}$} & $2^{-5}$ & 61504 & 10 & 0.12 & 3.72e-02 \\
& & $2^{-6}$ & 123008 & 10 & 0.20 & 1.85e-02 \\
& & $2^{-7}$ & 246016 & 8 & 0.30 & 9.24e-03 \\
& & $2^{-8}$ & 492032 & 8 & 0.50 & 4.61e-03 \\
& & $2^{-9}$ & 984064 & 8 & 0.78 & 2.31e-03 \\
& & $2^{-10}$ & 1968128 & 8 & 1.54 & 1.15e-03 \\
& & $2^{-11}$ & 3936256 & 8 & 3.18 & 5.76e-04 \\ \hline
\end{tabular}
\end{table}

\begin{table}[]
\caption{Results of GMRES-$\mathcal{P}$ for Example \ref{ex:example2} with various regularization parameters $\gamma$}
\label{tab:table2_RBD_precon_revised_compare}
\centering
\begin{tabular}{ccccccc}
\hline
\multirow{2}{*}{$\gamma$} & \multirow{2}{*}{$h$} & \multirow{2}{*}{$\tau$} & \multirow{2}{*}{DoF} & \multicolumn{3}{c}{GMRES-$\mathcal{P}$} \\ \cline{5-7}
& & & & Iter & CPU & $e_{h,\tau}$ \\ \hline
\multirow{7}{*}{$10^{-10}$} & \multirow{7}{*}{$2^{-5}$} & $2^{-5}$ & 61504 & 6 & 0.12 & 6.00e-06 \\
& & $2^{-6}$ & 123008 & 6 & 0.15 & 1.48e-06 \\
& & $2^{-7}$ & 246016 & 6 & 0.21 & 3.64e-07 \\
& & $2^{-8}$ & 492032 & 6 & 0.36 & 8.85e-08 \\
& & $2^{-9}$ & 984064 & 6 & 0.69 & 2.12e-08 \\
& & $2^{-10}$ & 1968128 & 8 & 1.73 & 5.19e-09 \\
& & $2^{-11}$ & 3936256 & 8 & 3.43 & 4.93e-09 \\ \hline
\multirow{7}{*}{$10^{-8}$} & \multirow{7}{*}{$2^{-5}$} & $2^{-5}$ & 61504 & 6 & 0.06 & 6.00e-06 \\
& & $2^{-6}$ & 123008 & 8 & 0.15 & 1.48e-06 \\
& & $2^{-7}$ & 246016 & 8 & 0.26 & 3.64e-07 \\
& & $2^{-8}$ & 492032 & 10 & 0.62 & 8.86e-08 \\
& & $2^{-9}$ & 984064 & 11 & 1.35 & 2.13e-08 \\
& & $2^{-10}$ & 1968128 & 13 & 3.09 & 5.24e-09 \\
& & $2^{-11}$ & 3936256 & 16 & 6.79 & 1.66e-09 \\ \hline
\multirow{7}{*}{$10^{-6}$} & \multirow{7}{*}{$2^{-5}$} & $2^{-5}$ & 61504 & 10 & 0.09 & 6.01e-06 \\
& & $2^{-6}$ & 123008 & 12 & 0.25 & 1.49e-06 \\
& & $2^{-7}$ & 246016 & 15 & 0.55 & 3.70e-07 \\
& & $2^{-8}$ & 492032 & 18 & 1.28 & 9.42e-08 \\
& & $2^{-9}$ & 984064 & 19 & 2.38 & 2.61e-08 \\
& & $2^{-10}$ & 1968128 & 17 & 4.08 & 8.56e-09 \\
& & $2^{-11}$ & 3936256 & 15 & 6.59 & 3.28e-09 \\ \hline
\multirow{7}{*}{$10^{-4}$} & \multirow{7}{*}{$2^{-5}$} & $2^{-5}$ & 61504 & 19 & 0.28 & 7.67e-06 \\
& & $2^{-6}$ & 123008 & 20 & 0.59 & 3.84e-06 \\
& & $2^{-7}$ & 246016 & 18 & 0.76 & 1.92e-06 \\
& & $2^{-8}$ & 492032 & 15 & 1.06 & 9.59e-07 \\
& & $2^{-9}$ & 984064 & 13 & 1.58 & 4.79e-07 \\
& & $2^{-10}$ & 1968128 & 11 & 2.68 & 2.40e-07 \\
& & $2^{-11}$ & 3936256 & 9 & 4.56 & 1.20e-07 \\ \hline
\multirow{7}{*}{$10^{-2}$} & \multirow{7}{*}{$2^{-5}$} & $2^{-5}$ & 61504 & 13 & 0.17 & 7.00e-04 \\
& & $2^{-6}$ & 123008 & 12 & 0.25 & 3.50e-04 \\
& & $2^{-7}$ & 246016 & 11 & 0.45 & 1.75e-04 \\
& & $2^{-8}$ & 492032 & 10 & 0.67 & 8.75e-05 \\
& & $2^{-9}$ & 984064 & 10 & 1.21 & 4.38e-05 \\
& & $2^{-10}$ & 1968128 & 8 & 2.01 & 2.19e-05 \\
& & $2^{-11}$ & 3936256 & 8 & 3.60 & 1.09e-05 \\ \hline
\multirow{7}{*}{$1$} & \multirow{7}{*}{$2^{-5}$} & $2^{-5}$ & 61504 & 10 & 0.12 & 3.72e-02 \\
& & $2^{-6}$ & 123008 & 10 & 0.23 & 1.85e-02 \\
& & $2^{-7}$ & 246016 & 10 & 0.40 & 9.24e-03 \\
& & $2^{-8}$ & 492032 & 10 & 0.72 & 4.61e-03 \\
& & $2^{-9}$ & 984064 & 10 & 1.28 & 2.31e-03 \\
& & $2^{-10}$ & 1968128 & 10 & 2.43 & 1.15e-03 \\
& & $2^{-11}$ & 3936256 & 10 & 4.80 & 5.76e-04 \\ \hline
\end{tabular}
\end{table}

\begin{figure}[h!]
    \centering
    \includegraphics[scale=0.75]{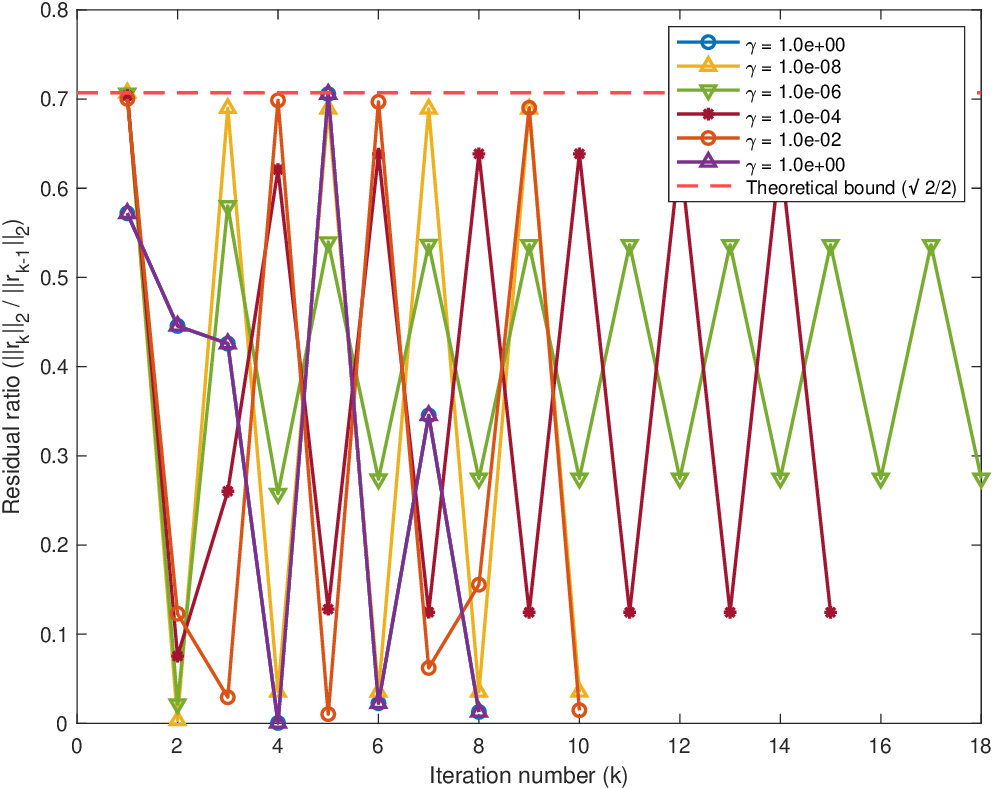}

        \caption{Convergence factor per iteration for Example~\ref{ex:example2} with various regularization parameters $\gamma$.}
        
    \label{fig:Example_2}
\end{figure}

					\end{example}

					\section{Conclusions}\label{sec:conclusion}
					In this study, we developed a novel PinT preconditioner for the all-at-once linear systems from optimal control problems with parabolic equations, using $\epsilon$-circulant matrices which can be efficiently implemented by fast Fourier transforms. Our approach enhances the convergence rate of the GMRES method, maintaining linearity with a theoretically estimated $\epsilon$, and remains independent of both the matrix size and the regularization parameter. Numerical experiments confirm the effectiveness and efficiency of our preconditioner. \correc{As part of our future work, we intend to generalize our preconditioning strategies to accommodate more complex problem settings—particularly those involving parabolic optimal control problems discretized in time using the second-order Crank–Nicolson scheme, as well as optimal control problems constrained by fractional time-dependent diffusive equations.}
					
					

					
					
					\section*{Acknowledgments}
					The work of Sean Y. Hon was supported in part by \correc{NSFC under grant 12401544,} the Hong Kong RGC under grant 22300921 and a start-up grant from the Croucher Foundation. The work of Xue-Lei Lin was partially supported by research grants: 12301480 from NSFC, HA45001143 from Harbin Institute of Technology, Shenzhen, HA11409084 from Shenzhen.

				\end{document}